\documentclass[11pt]{article}

\usepackage{amsmath, amsfonts, amssymb, amsthm}

\usepackage{natbib}

\usepackage{setspace}

\DeclareMathSymbol{\leqslant}{\mathalpha}{AMSa}{"36} % nicer `smaller or equal'

\DeclareMathSymbol{\geqslant}{\mathalpha}{AMSa}{"3E} % nicer `larger or equal'

\DeclareMathSymbol{\eset}{\mathalpha}{AMSb}{"3F}     % nicer `emptyset'

\renewcommand{\leq}{\;\leqslant\;}                   % redef. of < or =

\renewcommand{\geq}{\;\geqslant\;}                   % redef. of > or =

\newcommand{\ban}{\begin{align}}
\newcommand{\ean}{\end{align}}

\newcommand{\ba}{\begin{align*}}
\newcommand{\ea}{\end{align*}}

\newcommand{\be}{\begin{eqnarray*}}
\newcommand{\ee}{\end{eqnarray*}}

\newcommand{\ben}{\begin{eqnarray}}
\newcommand{\een}{\end{eqnarray}}

\theoremstyle{plain}
\newtheorem{theo}{Theorem}[section]

\newtheorem{lemma}[theo]{Lemma}
\newtheorem{propo}[theo]{Proposition}
\newtheorem{corollary}[theo]{Corollary}

\theoremstyle{definition}
\newtheorem{defi}[theo]{Definition}

\newtheorem{notation}[theo]{Notation}
\newtheorem{remark}[theo]{Remark}

\newtheorem{condition}[theo]{Condition}

\renewenvironment{proof}[1][] {{\bf Proof#1.} }{\hspace*{\fill}$\square$\medskip\par}

\begin{document}

\vglue20pt \centerline{\huge\bf Isotropic Ornstein-Uhlenbeck Flows }

\medskip

%\centerline{\huge\bf stochastic flows}

\bigskip

\bigskip

\centerline{by}

\bigskip

\medskip

\centerline{{\Large H.~van Bargen\footnotemark[1]  and  G.~Dimitroff\footnotemark[2]}}

\footnotetext[1]{Institut f\"ur Mathematik, MA 7-4, Technische Universit\"at Berlin, Stra\ss e des 17. Juni 136, D-10623 Berlin }

\footnotetext[2]{Fraunhofer ITWM, Fraunhofer-Platz 1, D-67663 Kaiserslautern }

\bigskip

\bigskip

{\leftskip=1truecm

\rightskip=1truecm

\baselineskip=15pt

\small

\noindent{\slshape\bfseries Summary.} Isotropic Brownian flows (IBFs) are a fairly natural class of stochastic flows which has  been studied extensively by various  authors. Their rich structure allows for explicit calculations in several situations and makes them a natural object  to start with if one wants to study  more general stochastic flows. Often the intuition gained by understanding the problem in the context of   IBFs transfers to more general situations. However, the obvious link between   stochastic flows,  random dynamical systems and ergodic theory  cannot be exploited in its full strength as the IBF does not have an invariant probability measure but rather an infinite one. Isotropic Ornstein-Uhlenbeck flows (IOUFs) are  in a sense localized IBFs and do have an invariant probability measure. The imposed linear drift destroys the translation invariance of the IBF, but many other important structure properties like the Markov property of the distance process remain valid and allow for explicit calculations in certain situations. The fact that IOUFs have invariant probability measures  allows one to apply techniques from random dynamical systems theory. We demonstrate this by applying the results of Ledrappier and Young to calculate the Hausdorff dimension of the statistical equilibrium of an IOUF. \bigskip


\bigskip


}

%%%%%% massblaetter

%\pagestyle{empty}
\newcommand{\oq}{{\langle}}
\newcommand{\ve}{{\varepsilon}}
\newcommand{\cq}{{\rangle}_t}
\newcommand{\dd}{{\mathrm{d}}}
\newcommand{\Dd}{{\mathrm{D}}}
\newcommand{\cd}{{\cdot}}
\newcommand{\nix}{{\varnothing}}
\newcommand{\N}{{\mathbb N}}
\newcommand{\Z}{{\mathbb Z}}
\newcommand{\R}{{\mathbb R}}
\newcommand{\Q}{{\mathbb Q}}
\newcommand{\E}{{\mathbb E}}
\renewcommand{\P}{{\mathbb P}}
\newcommand{\F}{{\cal F}}
\newcommand{\C}{{\mathbb C}}
\newcommand{\K}{{\mathbb K}}
\newcommand{\B}{{\cal B}}
\newcommand{\G}{{\cal G}}
\newcommand{\D}{{\cal D}}
\newcommand{\Zz}{{\cal Z}}
\newcommand{\Ll}{{\cal L}}
\newcommand{\A}{{\cal A}}
\newcommand{\Po}{{\cal P}}
\newcommand{\Sy}{{\cal S}}
\newcommand{\cZ}{{\cal Z}}
\newcommand{\M}{{\cal M}}
\newcommand{\Nn}{{\cal N}}
\newcommand{\p}{{\mathbf P}}
\newcommand{\X}{{\mathbb X}}
\newcommand{\interior}[1]{{\mathaccent 23 #1}}
\newcommand{\bem}{\begin{em}}
\newcommand{\eem}{\end{em}}
\def\eins{{\mathchoice {1\mskip-4mu\mathrm l}
{1\mskip-4mu\mathrm l}
{1\mskip-4.5mu\mathrm l} {1\mskip-5mu\mathrm l}}}
\newcommand{\olim}[1]{\begin{array}{c} ~\\[-1.4ex] \overline{\lim} \\[-1.35ex]
                       {\scriptstyle #1}\end{array}}
\newcommand{\plim}[2]{\begin{array}{c} #1\\[-1.175ex] \longrightarrow \\[-1.2ex]
                       {\scriptstyle #2}\end{array}}

\newcommand{\ulim}[1]{\begin{array}{c} ~\\[-1.175ex] \underline{\lim} \\[-1.2ex]
                       {\scriptstyle #1}\end{array}}

\newcommand{\prob}[2][]{\ensuremath{\mathbb{P}_{{#1}}\left[{#2}\right]}}

\renewcommand{\theequation}{\thesection.\arabic{equation}}

\section{Introduction}\label{intro}

Isotropic Brownian flows (IBFs) form a class of stochastic flows, which  is intimately connected to the Lebesgue measure in view of the following facts:
\begin{enumerate}
\item The one-point motion of an IBF has the Lebesgue measure as its invariant measure (modulo multiplicative constant).
\item The distribution of an IBF  is invariant under    rigid transformations   of $\R^d$, which also preserve, and in fact characterize, the Lebesgue measure in  $\R^d$ up to a multiplicative constant.
\end{enumerate}
The Gaussian measures are fairly  common  in probability theory   and therefore it is   natural to look at the family of stochastic flows which is  connected to the family of  centered  Gaussian measures on $\B(\R^d)$ in the above sense. The one-point motion of  such a flow should be an Ornstein-Uhlenbeck process, and its law must be invariant with respect to rotations of $\R^d$. The procedure of obtaining an Ornstein-Uhlenbeck process from a given Brownian motion, via an SDE with linear drift can also be performed  for flows.
Here is the definition of the isotropic Ornstein-Uhlenbeck flows, using this procedure:\\
Take  a complete probability space $\left(\Omega,\F,\P  \right) $ and let $F(t,x,\omega)$ be an isotropic Brownian field with $C^4$ covariance tensor $b$, i.e. $F$ generates an isotropic Brownian flow.

For the convenience of the reader we will first revise some important facts on Isotropic Brownian Flows and their generating fields.
The $\mathbb{R}^{d}$-valued random vector field $\left\{F(t,x): t\geq0, x\in\mathbb{R}^{d}\right\}$ is an isotropic Brownian field if
$(t,x)\mapsto F(t,x)$ is an almost surely continuous,  time homogeneous,  centered Gaussian process whose distribution is invariant with respect to rigid transformations (translation, reflection, rotation) of the spatial variable $x$. The covariance structure of $F$ is given by  $\mathit{cov}(F(s,x),F(t,y))=(s\wedge t)b(x-y)$, where $b:\R^d \to \R^{d\times d}$ is a so called isotropic covariance tensor. The isotropy of $F$ (distributional invariance with respect to rotations and reflections  in the spatial variable) implies that
 \begin{align}\label{isotropy_tensor}
 b(x)=O^Tb(Ox)O
 \end{align}
 for all $O$ from the orthogonal group $\mathcal{O}^d$. Following \cite{BaH86} we will throughout  assume  that $b$ has continuous and bounded derivatives of order up to $4$.  Some further details on the covariance tensor $b$ will be given in the beginning of Section~\ref{sec:genfacsec}.
 An IBF is defined to be the flow generated by the  Kunita-type SDE $\psi_{s,t}(x)=x+\int_s^tF(du,\psi_{s,u}(x))$. Let us note that the Proposition~\ref{loc_prop1} is valid for IBFs if one puts $c=0$ there and replaces the words \bem Ornstein-Uhlenbeck diffusion\eem \mbox{} by the words \bem Brownian Motion \eem\mbox{}.

For more details on  isotropic Brownian flows and their generating fields consult e.g.  \cite{BaH86} and \cite{LeJ85}.

 For a $c>0$ define the semimartingale field
$$V(t,x,\omega)=F(t,x,\omega)-\int\limits_0^t c x \dd s=F(t,x,\omega)- c xt \,.$$
It is an isotropic Brownian field with a linear drift towards the origin and its local characteristic $(b,l_c)$ (with $l_c(x):=cx$) belongs to the class $(B_{ub}^4,B_{ub}^\infty)$ (see \cite{Ku90} for a comprehensive study of semimartingales with spatial parameter). The obvious analogy with   the characterization  of an Ornstein-Uhlenbeck process as a Brownian motion in a quadratic potential motivates the following definition.
\begin{defi}\label{def_iouf}
An \bem  isotropic Ornstein-Uhlenbeck flow\eem, shortly  \bem IOUF \eem  is a stochastic flow of $C^{3,\delta}$ (for arbitrary $\delta>0$ ) diffeomorphisms $\phi$, generated by (via a Kunita-type SDE) a semimartingale field $V(t,x,\omega)$ as above, i.e.  $\phi_{s,t}(x)=x+\int_s^tV(du,\phi_{s,u}(x))$.
\end{defi}

Clearly, an IOUF is not translation invariant, but still retains most of the nice properties of  an IBF, like the most central fact that the distance process is a diffusion. In contrast to IBFs,   IOUFs have the very important feature that the one-point motion admits  an invariant \begin{em} probability \end{em} measure, and therefore one has at his disposal some powerful tools  from  random dynamical systems theory. To be more specific  the linearization of an IBF does not strictly fulfill the conditions of the multiplicative ergodic theorem as the flow does not have an invariant \begin{em} probability \end{em}  measure and therefore we only get some restricted statement instead ( the Lyapunov exponents are still well defined). On the other hand  the linearization of an IOUF satisfies the conditions of the multiplicative ergodic theorem and we have its statement in its full strength.  \\

To an IOUF as well as to an IBF one can associate a statistical equlibrium, which is the proper extension of the notion of an invariant measure to stochastic flows. It is a random measure on $\R^d$ (we consider $\R^d$-flows), which is stationary under the action  of the flow, and moreover the pullback of the invariant measure for the one point motion converges to the statistical equilibrium (more on this in Section~\ref{sec:dimsec}).  For an IBF we have the following: In the volume preserving case the statistical equilibrium is the Lebesgue measure itself, and in the case of a negative top Lyapunov exponent  Darling and LeJan (\cite{DaL88}) showed that  the statistical equilibrium is a zero measure. In the rest of the cases the statistical equlibrium is a non-trivial random measure being a.s.~singular to the Lebesgue measure. In the case of an IOUF the statistical equilibrium is a random probability measure being Dirac if the top Lyapunov exponent is strictly negative and a.s. diffuse if it is strictly positive.  Ledrappier and Young  were able to explicitly calculate (in terms of the Lyapunov exponents) the Hausdorff dimension of the statistical equilibrium for an RDS on a compact manifold  (\cite{LeY88}) and moreover, to link it to the entropy of the system, thus giving a precise meaning of the intuitively clear interplay between the entropy and dimensionality of the invariant measure of a dynamical system (\cite{LeY85_1}, \cite{LeY85_2}, \cite{LeY84}). In Section~\ref{sec:dimsec} we apply their result to IOUFs with the help of a compactification argument, which works precisely because an IOUF has an invariant probability measure, and consequently breaks down for isotropic Brownian flows. We conjecture that the Hausdorff dimension of an IBF with strictly positive top Lyapunov exponent can be obtained as a limit of the dimensions of the corresponding dimensions for the IOUFs with $c\to 0$, however we have not been able to prove this yet. The conjecture is trivially true for the extreme situations of an IBF with negative top Lyapunov exponent and a volume preserving IBF.

Another motivation for considering the IOUF is the natural question of letting an isotropic Brownian flow  evolve in a localizing potential,  i.e.~for some  field $U\colon \R^d \to \R$ consider  the IBF $\phi^U$ enforced by the potential $U$, that is, $\phi^U$ is generated via an SDE driven by the semimartingale field
\begin{align*}
V^U(t,x,\omega)=F(t,x,\omega)-\nabla U(x)t\,.
\end{align*}
The IOUFs correspond to a quadratic potential $U=\frac{c}{2}|x|^2$.
However, adding a drift  typically will destroy many of the  nice properties of an IBF.
Clearly, the IOUFs are not translation invariant, but still  retain most of the rich structure  of the IBFs, like as already mentioned, the most central fact that the distance process is a diffusion.
Moreover the IOUFs can be considered as localized IBFs as  their one-point motions have an invariant \begin{em} probability \end{em} measures. As we shall see, they are also localized as flows since they posses weak  random attractors.  \\
The paper is organized as follows: In the next section we provide some  ``infrastructure'' facts on IOUFs.  In Section 3 we introduce the notion of statistical equilibrium and give its Hausdorff dimension for an IOUF. In Section 4 we give a brief introduction to the notion of random attractors and then establish the existence of  weak attractors for IOUFs.

\section{General facts on IOUF}
\label{sec:genfacsec}

Let $\left(\Omega,\F,(\F_s^t)_{0\le s\le t},\P \right)$ be a filtered probability space satisfying the usual conditions i.e. every $\F_s^t$ is complete and the filtration is right-continuous in  $s$ and $t$. Let $\phi$ be an IOUF, generated as above by the field $F(t,x,\omega)- c xt$.
\begin{notation} We will write  $\phi_t$ for $\phi_{0,t}$.\end{notation}
 The distributions of $F$ and $\phi$ are uniquely determined by the \bem covariance tensor  \eem
\begin{align*}
b:\R^d\times\R^d\to \R^{d\times d}  \, , \,\,\, (x,y)\mapsto \E F(1,x)F(1,y)\textnormal{ and }c>0.
\end{align*}
According to Yaglom \cite[Section 4]{Yag57} (see also
\cite{BaH86}) the isotropy property (\ref{isotropy_tensor}) implies that  $b$   has necessarily the form
\begin{align}\label{pre_ibf2}
b_{i,j}(x)= \left\{ \begin{array}{ll} \left(B_L(|x|)-B_N(|x|)\right)\frac{x_ix_j}{|x|^2}+\delta_{i,j}B_N(|x|) & \text{for }x\ne 0  \\ [0.4cm]
\delta_{i,j}B_L(0)=\delta_{i,j}B_N(0) & \text{for }x= 0\,,
\end{array} \right.
\end{align}
where $B_L$ and $B_N$ are the so called longitudinal and  transversal correlation  functions defined by
\begin{align*}
&B_L(r)=b_{p,p}(re_p), \,\,\, r\ge 0 \,\,\,\text{ and }\\
&B_N(r)=b_{p,p}(re_q), \,\,\, r\ge 0,\,\,\, p\ne q \,.
\end{align*}
In the above formulas $ e_i $ denotes the $ i $-th standard basis vector in $ \R^d $.  In particular, (\ref{isotropy_tensor}) implies $b(0)=k\cdot \text{Id}_{\R^d}$ for some positive constant $k$. Without loss of generality, as it is done in~\cite{BaH86},  we  assume that $b(0)= \text{Id}_{\R^d}$. In this case the  one-point-motion is a standard Brownian motion and we also get $B_L(0)=B_N(0)=1$. The general case  $b(0)=k\cdot \text{Id}_{\R^d}$ can be obtained via rescaling the time by the  constant factor $\sqrt{k}$.
$B_L$ and $B_N$ are bounded $C^4$ functions with bounded derivatives up to order four. Set
\begin{align*}
\beta_L&:=-\partial_p\partial_pb_{p,p}(0)=-B_L''(0)\\
\beta_N&:=-\partial_q\partial_q b_{p,p}(0)=-B_N''(0) \,\,\, q\ne p\,,
\end{align*}
where $B_L''(0)$ and $B_N''(0)$ denote the right-hand second derivatives of $B_L$ and $B_N$ at zero. We refer the reader to  \cite{BaH86} and  \cite{LeJ85} for more details.
\\[0.4cm]
\begin{propo}\label{loc_prop1}
Let $\phi$ be an isotropic Ornstein-Uhlenbeck flow. Then for $x,y\in\R^d$
\begin{enumerate}
\item[\textnormal{(a)}] $\phi$ is Brownian and  isotropic (invariant with respect to orthogonal transformations).
\item[\textnormal{(b)}] The one-point motion $\{\phi_t(x)\colon t\in\R\} $ of $\phi$ is an Ornstein-Uhlenbeck diffusion with generator
\begin{align*}
\mathcal{L}:=-c\sum_{i=1}^d x_i\frac{\partial}{\partial x_i}+\frac{1}{2}\sum_{i=1}^d\frac{\partial^2}{\partial x_i^2} =-c\sum_{i=1}^d x_i\frac{\partial}{\partial x_i}+\frac{1}{2}\triangle\,.
\end{align*}
\item[\textnormal{(c)}] The difference process $\{\phi_t(x)-\phi_t(y) : t \in \R_+\}$ is a diffusion with generator
\begin{align*}
\mathcal{L}_d:=-c\sum_{i=1}^d x_i\frac{\partial}{\partial x_i}+\sum_{i,j=1}^d (\delta_{ij}-b_{i,j}(x)) \frac{\partial^2}{\partial x_i\partial x_j}\,.
\end{align*}
\item[\textnormal{(d)}] The distance process $\{|\phi_t(x)-\phi_t(y)| : t \in \R_+\}$ is a diffusion with generator
\begin{align*}
\mathcal A=(1-B_L(r))\frac{\dd^2}{\dd r^2}+\big ( (d-1)\frac{1-B_N(r)}{r}-c r \big )\frac{\dd}{\dd r}\,.
\end{align*}
\end{enumerate}
\end{propo}
\begin{proof}
\textbf{(a)} Recall, that a stochastic flow is called \bem Brownian \eem  if it has independent increments (see \cite{Ku90}, Section 4.1, page 116), i.e.~for all $n \in \N$ and $0\le t_1\le t_2\le \dots\le t_n$ the  random variables $(\phi_{t_1,t_2}, \dots, \phi_{t_{n-1},t_n})$  are independent. This property is obvious as the field $F$ has independent increments.    \\
The  invariance with respect to orthogonal transformations  is a direct consequence of the  invariance of $F(t,x) $ and the drift with respect to these operations.
\textbf{(b)} Using Definition~\ref{def_iouf} we get for the joint quadratic variation of the components of the one-point motion
\begin{align*}
\langle \phi_{\cdot}^i(x),\phi_{\cdot}^j(x)\rangle_t=\langle \int\limits_0^\cdot F^i(\dd s,\phi_{s}(x)),\int\limits_0^\cdot F^j(\dd s,\phi_{s}(x))\rangle_t=b_{i,j}(0)t=\delta_{ij}t\,.
\end{align*}
Let $f\in C^2_0(\R^d:\R)$ be arbitrary. An application of  It\^o's formula gives
\begin{align*}
&f(\phi_t(x))-f(x)=\sum_{i=1}^d\int\limits_0^t \partial_i f(\phi_s(x))\dd\phi_{s}^i(x)+\frac{1}{2} \sum_{i,j=1}^d\int\limits_0^t \partial^2_{ij} f(\phi_s(x))\dd \langle\phi_{\cdot}^i(x), \phi_{\cdot}^j(x)\rangle_s\\[0.4cm]
&=\sum_{i=1}^d\int\limits_0^t \partial_i f(\phi_s(x)) F^i(\dd s,\phi_{s}(x)) - c\sum_{i=1}^d\int\limits_0^t \phi_s^i(x)\partial_i f(\phi_s(x))\dd s +\frac{1}{2}\sum_{i,j=1}^d \delta_{ij}\int\limits_0^t \partial^2_{ij} f(\phi_s(x))\dd s\,,
\end{align*}
that is
$t\mapsto f(\phi_t(x))-f(\phi_0(x))-\int\limits_0^t (\mathcal{L}f)(\phi_s(x))\dd s$
is a continuous martingale for all $x\in \R^d$, and therefore $\phi_t(x)$ is a diffusion with generator $\mathcal{L}$ (see \cite{KaShr88}, Propositions 4.6 and 4.11,  Section 5.4.B). \\[0.5cm]
\textbf{(c)} Let $x,\,y \in \R^d$, $x\ne y$ be arbitrary but fixed. We clearly have  the equality:
\begin{align*}
& \phi_{t}(x)-\phi_{t}(y)=x-y -c  \int\limits_0^t(\phi_{t}(x)-\phi_{t}(y))\dd s+\int\limits_0^t F(\dd s,\phi_s(x))-\int\limits_0^t F(\dd s,\phi_s(y))\,.
\end{align*}
Observe the following computation for   the joint quadratic variation of the difference process
\begin{align*}
&\langle \phi_{\cdot}^i(x)-\phi_{\cdot}^i(y)\,,\,\phi_{\cdot}^j(x)-\phi_{\cdot}^j(y)\rangle_t=\langle \phi_{\cdot}^i(x)\,,\,\phi_{\cdot}^j(x)\rangle_t+\langle \phi_{\cdot}^i(y)\,,\,\phi_{\cdot}^j(y) \rangle_t-2\langle \phi_{\cdot}^i(x)\,,\,\phi_{\cdot}^j(y)\rangle_t\\[0.4cm]
&=2 \int\limits_0^t(\delta_{ij}-b_{i,j}(\phi_s(x)-\phi_s(y)))ds\,.
\end{align*}
For arbitrary $f \in C^2_0(\R^d:R)$ with the help of   It\^o's formula we obtain
\begin{align*}
&f(\phi_t(x)-\phi_t(y))-f(x-y)=\sum_{i=1}^d\int\limits_0^t \partial_i f(\phi_s(x)-\phi_s(y))\dd(\phi_{s}^i(x)-\phi_s^i(y))\\[0.4cm]
&+ \frac{1}{2}\sum_{i,j=1}^d\int\limits_0^t \partial^2_{ij} f(\phi_s(x)-\phi_s(y))\dd \langle\phi_{\cdot}^i(x)-\phi_{\cdot}^i(y), \phi_{\cdot}^j(x)-\phi_{\cdot}^j(y)\rangle_s\\[0.4cm]
&=\sum_{i=1}^d\int\limits_0^t \partial_i f(\phi_s(x)-\phi_s(y))F^i(\dd s,\phi_s(x))-\sum_{i=1}^d\int\limits_0^t \partial_i f(\phi_s(x)-\phi_s(y))F^i(\dd s,\phi_s(y))\\[0.4cm]
&-c\sum_{i=1}^d(\phi_{t}^i(x)-\phi_{t}^i(y))\partial_i f(\phi_s(x)-\phi_s(y))\dd s\\[0.4cm]
&+\frac{1}{2}\sum_{i,j=1}^d\int\limits_0^t 2(\delta_{ij}-b_{i,j}(\phi_s(x)-\phi_s(y)))\partial^2_{ij} f(\phi_s(x)-\phi_s(y))\dd s \,,
\end{align*}
that is
$t\mapsto f(\phi_t(x)-\phi_t(y))-f(x-y)-\int\limits_0^t (\mathcal{L}_df)(\phi_{s}(x)-\phi_s(y))\dd s$ is a continuous martingale for all $x\in\R^d$ and we can conclude as before.
\\
\textbf{(d)}Let $x\ne y \in \R^d$  be arbitrary but fixed.\\
Denote by $r_t$ the distance process $|\phi_t(x)-\phi_t(y)|$. The homeomorphic property of $\phi$ implies that $\P(r_t>0\,, \textrm{for all } t \in \R_+)=1$. Set  $q_t:=\phi_t(x)-\phi_t(y)$. For arbitrary  $f\in C^2_0(\R :\R)$  the   It\^o formula implies that
\begin{align*}
&f(r_t)-f(r_0)=\underbrace{\sum_{i=1}^d\int\limits_0^t f'(r_s)\frac{q^i_s}{r_s}F^i(\dd s,\phi_s(x))-\sum_{i=1}^d\int\limits_0^t f'(r_s)\frac{q^i_s}{r_s}F^i(\dd s,\phi_s(y))}_{M_t}\\[0.4cm]
&+ \int\limits_0^t\Big[ ( -c r_s+\frac{d-1}{r_s}-\sum_{i=1}^d\frac{b_{i,i}(q_s)}{r_s}+\sum_{i,j=1}^d\frac{q_s^iq_s^j}{r_s^3} b_{i,j}(q_s) )f'(r_s)+ (1-\sum_{i,j=1}^d b_{i,j}(q_s)\frac{q_s^iq_s^j}{r_s^2})f''(r_s)\Big]\dd s\,.
\end{align*}

Using the equalities
\begin{align*}
\sum_{i,j=1}^d b_{i,j}(x)x^ix^j&=(B_L(|x|)-B_N(|x|))\sum_{i,j=1}^d \frac{(x^ix^j)^2}{|x|^2} +\sum_{i,j=1}^d\delta_{ij}B_N(|x|)x^ix^j\\[0.4cm]
&=(B_L(|x|)-B_N(|x|))|x|^2+B_N(|x|)|x|^2=B_L(|x|)|x|^2
\end{align*}
and
\begin{align*}
\sum_{i=1}^d b_{i,i}(x)&=(B_L(|x|)-B_N(|x|))\sum_{i=1}^d \frac{(x^i)^2}{|x|^2} +\sum_{i=1}^d\delta_{ii}B_N(|x|)=B_L(|x|)+(d-1)B_N(|x|)\,
\end{align*}
we obtain
\begin{align*}
&f(r_t)-f(r_0)=M_t+ \int\limits_0^t\Big[( (d-1)\frac{1-B_N(r_s)}{r_s}-c r_s )f'(r_s)+(1-B_L(r_s))f''(r_s) \Big]\dd s\,,
\end{align*}
that is
\begin{align*}
f(r_t)-f(r_0)-\int_0^t (\mathcal{A}f)(r_s)\dd s
\end{align*}
is a continuous martingale and  standard results about the martingale problem imply that $(r_t)_{t\ge 0}$ is a diffusion with generator $\mathcal{A}$.
\end{proof}
The following corollary will be used throughout the paper. It has been proven in a very general setting by Cranston, Scheutzow and Steinsaltz in \cite{CSS00}, however they impose too strong boundedness assumptions on the local characteristic of the generating field $V$ to be directly applicable to the case of an IOUF. In particular, they assume that the drift is uniformly bounded which is obviously not fulfilled in our case.    Their proof however does work under much milder conditions which include the case of an IOUF. Therefore, their statement formally does not cover our case, however their proof does. Recent joint work of one of the authors and Michael Scheutzow contains the following bound under milder conditions  making it directly applicable to our case. For the sake of completeness we give the proof for the case of an IOUF.
\begin{corollary}\label{iouf_gen_cor_0}
Let $\phi$ be an IOUF and let $x,y\in \R^d$, $x\ne y$ be arbitrary. Then there is a standard Brownian motion $\left\{ B_{t} : t\ge 0\right \}$ (depending on $x$ and $y$!) and constants $\sigma,\,\lambda>0$  such that $\P$-almost surely
$$|\phi_t(x)-\phi_t(y)|\le |x-y|\exp\left\{ \sigma B_t^\star + \lambda t  \right\}\,\,\, \text{ for all } t\ge0\,,$$
where as usual $B^\star_t:=\sup_{s\le t}B_s$.
\end{corollary}
\begin{proof}
According to Proposition \ref{loc_prop1} (d) we can view  $\left\{r_t:=|\phi_t(x)-\phi_t(y)|: t\ge 0\right\}$ as the unique strong solution of
\begin{align*}%\label{iouf_gen_cor_01}
r_t=|x-y|+\int\limits_0^t\sqrt{2(1-B_L(r_s))}\,\dd W_s+\int\limits_0^t\left( (d-1)\frac{1-B_N(r_s)}{r_s}-cr_s\right)\,\dd s\,.
\end{align*}
where $\left\{ W_{t} : t\ge 0\right\}$ is a standard Brownian motion  (depending on $x$ and $y$!).
The asymptotic expansions of $B_L$ and $B_N$ around  zero (see \cite{BaH86}) provide constants $a,b>0$, such that
\begin{align}\label{iouf_gen_cor_02}
\sup\limits_{u\ge 0}\frac{1-B_N(u)}{u^2}\le a \,\,\,\text{ and }\,\,\, \sup\limits_{u\ge 0}\frac{1-B_L(u)}{u^2}\le b\,.
\end{align}
With
$$M_t=\int\limits_0^t\sqrt{\frac{2(1-B_L(r_s))}{r_s^2}}\dd W_s$$
the SDE for $r_t$ can be rewritten as
\begin{align*}%\label{iouf_gen_cor_04}
r_t=|x-y|+\int\limits_0^t r_s\,\dd M_s+\int\limits_0^t\Big( (d-1)\frac{1-B_N(r_s)}{r_s^2}-c\Big)r_s\,\dd s\,.
\end{align*}
Consider the process $(Y_t)$, solving the SDE
\begin{align*}
Y_t=|x-y|+\int\limits_0^t Y_s\,\dd M_s+\int\limits_0^t \left((d-1)a-c\right)Y_s \,\dd s\,,
\end{align*}
that is
\begin{align*}%\label{iouf_gen_cor_05}
Y_t=|x-y|\exp\Big\{ M_t+ t \left((d-1)a-c\right)-\frac{1}{2}\langle M \rangle_t\Big \}\,.
\end{align*}
An application of the comparison theorem as stated in \cite{Pr04} (Theorem 54 on page 324, Chapter V, Section 9)  implies
\begin{align}\label{iouf_gen_cor_06}
\P\left( r_t \le Y_t \text{ for all } t\ge 0\right) =1\,,
\end{align}
as  we have $ (d-1)\frac{1-B_N(u)}{u^2}-c \le (d-1)a-c\,.$
According to (\ref{iouf_gen_cor_02}) for the quadratic variation of $M$ we have $\langle M \rangle_t\le 2bt\,.$
Recall that there is a standard Brownian motion $\left( B_t \right)_{t\ge 0}$ such that \\ $\P(M_t=B_{\langle M \rangle_t}\text{ for all }t<\infty)=1$. Thus for $Y_t$ we get $\P$-almost surely:
\begin{align*}
Y_t\le |x-y|\exp\left[ B_{\langle M \rangle_t}+ t \left((d-1)a-c\right)\right]\le |x-y|\exp\left[ B_{\langle M \rangle_t}^\star+ t \left((d-1)a-c\right)\right] \\%[0.4cm]
\le  |x-y|\exp\left[ B_{2bt}^\star+ t \left((d-1)a-c\right)\right]= |x-y|\exp\left[ \sqrt{2b}\sup_{s\le t}\tilde{B}_{s}+ t \left((d-1)a-c\right)\right]\,,
\end{align*}
for the Brownian motion $\tilde{B}_s=(2s)^{-\frac{1}{2}}B_{2bt}$, In view of  (\ref{iouf_gen_cor_06}) we are done
\end{proof}
In the next proposition we give the Lyapunov spectrum of an isotropic Ornstein-Uhlenbeck flow. Actually strictly speaking we are allowed to speak about Lyapunov exponents only  for random dynamical systems (RDS). However to a stochastic flow $\phi$ with stationary increments generated via a Kunita type SDE and satisfying some regularity conditions one can canonically construct a corresponding two-sided RDS $\varphi$ having the same distribution.   The proof can be found in \cite{AS95} without the extension from one sided time to $\R$, which is however canonical. %The complete proof is also quoted in (\cite{Dim06} Proposition 2.2.1.).
We make this  precise for an IOUF in the following remark.
\begin{remark}\label{sf_rds} Let $\phi$ be an IOUF. There is a perfect cocycle $\varphi$ on $\R^d$ with time $\R$ over some metric dynamical system  $(\tilde \Omega,\tilde \F,( \theta_t)_{t\in \R},\tilde \P)$ such that the distributions of $\big\{\phi_{s,t} : s,t\ge 0\big\}$ and $\big\{\varphi(t-s, \theta_s\tilde \omega) : s,t\ge 0\big\}$ coincide.\\
\end{remark}
From now on we shall identify the flow with the associated RDS.
\begin{propo}\label{loc_prop2}
The isotropic Ornstein-Uhlenbeck flow $\phi$ in $\R^d$ has $d$ different Lyapunov exponents, which are given by
\begin{align*}
\lambda_i:=(d-i)\frac{\beta_N}{2}-i\frac{\beta_L}{2}-c.
\end{align*}
In particular they all have simple multiplicities.
\end{propo}

\begin{proof}
Consider the isotropic Brownian flow $\psi$, generated by the field $F$ via the SDE
\begin{align*}
\psi_{s,t}(x)=x+\int\limits_s^t F(\dd u,\psi_{s,u}(x)) \text{ for all } 0\le s\le t,\, x\in \R^d
\end{align*}
and the IOUF $\phi$:
\begin{align*}
\phi_{s,t}(x)=x+\int\limits_s^t F(\dd u,\phi_{s,u}(x))-c\int\limits_s^t \phi_{s,u}(x) \dd u\,.
\end{align*}
Let $\Dd_x F(t,x)$ denote the Jacobi matrix of $F(t,x)$. Recall the convention $\phi_t:=\phi_{0,t}$ and $\psi_t:=\psi_{0,t}$. According to  Theorem 3.3.4 from \cite{Ku90} we have
\begin{align*}
\Dd_x\psi_{t}(x)=\text{Id}_{\R^d}+\int\limits_0^t \Dd_x F(\dd u,\psi_{u}(x))\Dd_x\psi_{u}(x)
\end{align*}
and
\begin{align*}
\Dd_x\phi_{t}(x)=\text{Id}_{\R^d}+\int\limits_0^t \Dd_x F(\dd u,\phi_{u}(x))\Dd_x\phi_{u}(x)-c\int\limits_0^t \Dd_x\phi_{u}(x)\dd u\,.
\end{align*}
Set $\gamma_t(x):=e^{-ct}\Dd_x\psi_t(x)$. It\^o's formula implies  that $\gamma$ solves the SDE
\begin{align*}
\gamma_t(x)=\text{Id}_{\R^d}+\int\limits_0^t \Dd_x F(\dd s,\psi_s(x))\gamma_s(x)-c\int\limits_0^t \gamma_s(x)\dd s\,.
\end{align*}
For arbitrary but fixed $x\in \R^d$ consider the martingale fields
\begin{align*}
&M^x(t,y)=\int\limits_0^t \Dd_x F(\dd s,\phi_s(x))y\,\colon \R_+\times \R^{d\times d} \times \Omega \to \R^{d\times d}\text{ and}\\
&N^x(t,y)=\int\limits_0^t \Dd_x F(\dd s,\psi_s(x))y\,\colon \R_+\times \R^{d\times d} \times \Omega \to \R^{d\times d}\,.
\end{align*}
Straightforward calculation of their quadratic variations shows that
\begin{align*}
\langle M^x_{ij}(\cdot,y),M^x_{pq}(\cdot,z) \rangle_t= \langle N^x_{ij}(\cdot,y)\,,\,N^x_{pq}(\cdot,z) \rangle_t=t\sum_{k,l=1}^d \partial_{k}\partial_{l}b_{i,p}(0)y_{kj}y_{lq}
\end{align*}
and therefore   $N^x$ and $M^x$ have the same Gaussian distribution. Clearly, the equations for $e^{-ct}D_x\psi_t(x)$ and $D_x\phi_t(x)$ can be written with the help of the fields $M^x$ and $N^x$:
\begin{align}\label{iouf_gen_prop2.1}
\begin{split}
\gamma_t(x)&=\text{Id}_{\R^d}+\int\limits_0^t M^x(\dd s,\gamma_s(x))-c\int\limits_0^t \gamma_s(x)\dd s\,, \\
\Dd_x\phi_t(x)&=\text{Id}_{\R^d}+\int\limits_0^t N^x (\dd s,\Dd_x\phi_s(x))-c\int\limits_0^t \Dd_x\phi_s(x)\dd s\,
\end{split}
\end{align}
and thus  $\big\{\gamma_t(x) : t\ge 0\big\}$ and $\big\{\Dd_x\phi_t(x) : t\ge 0\big \}$ have the same distributions for all $x \in \R^d$.\\
The  Lyapunov exponents of $\psi$ have simple multiplicities and are given by  $(d-i)\frac{\beta_N}{2}-i\frac{\beta_L}{2}$ (see \cite{BaH86}) and therefore almost surely
\begin{align*}
\lim_{t \to +\infty} \frac{1}{2t}\log\left[\sigma_i\left(\Dd_x\psi_t(x)\right)\right]=(d-i)\frac{\beta_N}{2}-i\frac{\beta_L}{2}\,,
\end{align*}
where $\sigma_i(M)$ denotes the $i$-th characteristic value of the matrix $M$ (the $i$-th eigenvalue of $M^TM$). According to the multiplicative ergodic theorem the corresponding a.s.~limits for $\phi$ also exist (see  \cite{MoSch99}) and satisfy
\begin{align*}
\lambda_i&:=\lim_{t \to +\infty} \frac{1}{2t}\log(\sigma_i(\Dd_x\phi_t(x)))=\lim_{t \to +\infty} \frac{1}{2t}\log(\sigma_i(\gamma_t(x)))=\lim_{t \to +\infty} \frac{1}{2t}\log(\sigma_i(\Dd_x\psi_t(x))))-c\\[0.4cm]
&=(d-i)\frac{\beta_N}{2}-i\frac{\beta_L}{2}-c\,,
\end{align*}
\end{proof}
In the last proposition in this section we give the transience/recurrence modes of the $\R_+$-valued diffusion $r_t(x,y):=|\phi_t(x)-\phi_t(y)|$. The strong inwards drift clearly implies \\$\P(\lim\limits_{t \to \infty}|\phi_t(x)-\phi_t(y)| =+\infty)=0$ and therefore ``transience'' means  almost sure convergence to zero, i.e. $\P(\lim\limits_{t \to \infty}r_t(x,y) =0)=1$.
\begin{propo}\label{loc_prop3}
Let $\lambda_1$ denote the top Lyapunov exponent of the IOUF $\phi$. Then for arbitrary $x,\,y\in \R^d$ the   diffusion $r_t=r_t(x,y)$ on $(0,+\infty)$ is
\begin{enumerate}
\item[\textnormal{(a)}] recurrent if $\lambda_1\ge 0$, i.e.
$$\P(\liminf_{t\to \infty}r_t(x,y)=0)=\P(\limsup_{t\to \infty}r_t(x,y)=+\infty)=1$$ %and
%$$\E\big[ \inf\{t> 0 : r_t(x,y)=|x-y|\} \big]<\infty\,.$$
%\item~ null recurrent if $\lambda=0$, i.e.
%$$\P(\liminf_{t\to \infty}r_t(x,y)=0)=\P(\limsup_{t\to \infty}r_t(x,y)=+\infty)=1$$ and
%$$\E\big[ \inf\{t> 0 : r_t(x,y)=|x-y|\} \big]=\infty\,.$$
\item[\textnormal{(b)}] transient if $\lambda_1<0$, that is
$$\P(\lim_{t\to \infty}r_t(x,y)=0)=1\,.$$
\end{enumerate}
Further, the \bem speed measure \eem of the diffusion $ r_t $ is finite if and only if $\lambda>0$. In this case the unique invariant probability measure of $ r_t $ is given by
$$\frac{1}{\int_0^\infty m(u)\dd u} m(x)\dd x\,,\,\,\,\text{ where }\,\,\,\,m(x)= \frac{1}{1-B_L(x)}\exp \Big [ \int\limits_{1}^x \frac{(d-1)(1-B_N(y))-cy^2}{y(1-B_L(y))}\dd y\Big ]\,.$$
\end{propo}
\begin{proof}
Observe that the exit time $S:=\inf\{t\ge 0 : r_t(x,y)=0 \text{ or } r_t(x,y)=+\infty\}$ is infinite almost surely. Indeed, for $ x\ne y $, $\P(r_t(x,y)>0 \text{ for all } t\ge 0)=1$ by the homeomorphic property of the flow $\phi$. Clearly  $\P(r_t(x,y)<+\infty \text{ for all } t\ge 0)=1$. % is a direct consequence of the fact that the coefficients in the SDE defining $r_t$ are globally Lipschitz and therefore the solution is globally defined and pathwise unique (see e.g. \cite{KaShr88}, Section5.2.B, Theorem 2.9).
 The density  $m(x)$  of the \bem speed measure \eem (see \cite{KaShr88}, Section 5.5.C, page 343) of $r_t(x,y)$    with respect to the Lebesgue measure satisfies the equation  $\mathcal{A}^\ast m=0$ and  can be calculated explicitly (see e.g. \cite{KaShr88}, Section 5.5.C, eq. 5.51). Up to a multiplicative constant it is given by
\begin{align*}
m(x)= \frac{1}{1-B_L(x)}\exp \Big [ \int\limits_{1}^x \frac{(d-1)(1-B_N(y))-cy^2}{y(1-B_L(y))}\dd y\Big ]\,.
\end{align*}
The \bem scale function \eem (see \cite{KaShr88}, Section 5.5.C, page 343) $s(x)$ is a solution of $\mathcal{A}u=0$ and  up to  multiplicative and additive constants is given by (see e.g. \cite{KaShr88}, Section 5.5.B, eq. 5.42)
\begin{align*}
s(x)=\int\limits_{1}^x \text{exp}\Big[  -\int\limits_{1}^y \frac{(d-1)(1-B_N(z))-cz^2}{z(1-B_L(z))} \dd z\Big]\dd y \,.
\end{align*}
Observe that since there is a $\delta \in (0,1)$, such that $-1+\delta<B_L(z),B_N(z)<1-\delta$ for all $z\ge 1$ (see \cite{BaH86}, Remark 2.18) we have for all $z\ge 1$
 \begin{align}\label{iouf_gen_prop3.1}
&- \frac{(d-1)(1-B_N(z))-cz^2}{z(1-B_L(z))}= - \frac{(d-1)(1-B_N(z))}{z(1-B_L(z))}+ \frac{cz}{1-B_L(z)}\ge -\frac{2(d-1)}{\delta z}+\frac{c}{2}z .
\end{align}
Therefore for $x\ge 1$
\begin{align}\label{loc_uminf}
s(x)  \ge \int\limits_{1}^x \text{exp}\Big[ \int\limits_{1}^y   \frac{c}{2}z-\frac{2(d-1)}{\delta z}\, \dd z\Big]\dd y=\int\limits_{1}^x \text{exp}\Big[  \frac{c}{4}y^2-\frac{2(d-1)}{\delta}\ln(y)-\frac{c}{4} \Big] \dd y \overset{x \to +\infty } \longrightarrow +\infty  \,.
\end{align}
In order to obtain the asymptotics of $s(x)$ as $x\to +0$, recall that for $r \to 0$ we have (see \cite{BaH86})
\begin{align*}
B_L(r)=1-\frac{\beta_L}{2}r^2+O(r^4)\,\,\text{ and }\,\,  B_N(r)=1-\frac{\beta_N}{2}r^2+O(r^4) \,
\end{align*}
and therefore
\begin{align*}
- \frac{(d-1)(1-B_N(z))-cz^2}{z(1-B_L(z))}= -\frac{1}{z}\frac{(d-1)\frac{\beta_N}{2}-c+O(z^2)}{\frac{\beta_L}{2}+O(z^2)}\,.
\end{align*}
Since
\begin{align*}
\frac{(d-1)\frac{\beta_N}{2}-c+O(z^2)}{\frac{\beta_L}{2}+O(z^2)}=\frac{(d-1)\frac{\beta_N}{2}-c}{\frac{\beta_L}{2}}+O(z^2)
\end{align*}
 we have that for  $x\to 0$,
\begin{align*}
s'(x)\sim  \text{exp}\Big [  -\frac{(d-1)\frac{\beta_N}{2}-c}{\frac{\beta_L}{2}} \int\limits_{1}^x \frac{1}{z} \dd z \Big]= x^ {-\frac{(d-1)\frac{\beta_N}{2}-c}{\frac{\beta_L}{2}}} \,,
\end{align*}
which is integrable around zero if and only if
\begin{align}\label{loc_umzero}
\frac{(d-1)\frac{\beta_N}{2}-c}{\frac{\beta_L}{2}}<1 \Longleftrightarrow \lambda_1<0 \,.
\end{align}
We now consider the density of the speed measure $m(x)$. Similarly as above, for all $x\ge 1$ we have  $-1+\delta<B_L(x),B_N(x)<1-\delta$ and therefore
\begin{align*}
m(x)&= \frac{1}{1-B_L(x)}\exp \Big [ \int\limits_{1}^x \frac{(d-1)(1-B_N(y))-cy^2}{y(1-B_L(y))}\dd y \Big ] \le \frac{1}{\delta}\exp \Big [ \int\limits_{1}^x \frac{2(d-1)}{\delta}\frac{1}{y}-\frac{c}{2}y\,\dd y \Big ]\\[0.4cm]
&=\frac{1}{\delta}\exp \Big [ \frac{2(d-1)}{\delta} \ln(x)- \frac{c}{4}x^2+\frac{c}{4} \Big ] \,,
\end{align*}
which is certainly integrable around $+\infty$.\\
For $x\to +0$
 \begin{align*}
m(x)\sim \frac{1}{\frac{\beta_L}{2}x^2}\exp \Big [  \frac{(d-1)\frac{\beta_N}{2}-c}{\frac{\beta_L}{2}}\int\limits_{1}^x \frac{1}{y}\dd y \Big ] = \frac{2}{\beta_L} x^{\frac{(d-1)\frac{\beta_N}{2}-c}{\frac{\beta_L}{2}}-2}\,,
\end{align*}
which is integrable around zero if and only if
$$\frac{(d-1)\frac{\beta_N}{2}-c}{\frac{\beta_L}{2}}-2>-1 \Longleftrightarrow \frac{(d-1)\frac{\beta_N}{2}-c}{\frac{\beta_L}{2}}>1 \Leftrightarrow \lambda_1>0\,.$$
In the case $\lambda_1=0$ the speed measure puts infinite charge on any  interval of the form $[0,\epsilon]$.
Combining these results about the behavior of the scale function near zero and infinity we get (see \cite{KaShr88}, Proposition 5.5.22)
\begin{enumerate}
\item (Transience) If $\lambda_1 < 0$ then $s(+\infty)=+\infty$ and $s(+0)>-\infty$, and therefore
\begin{align*}
\P \big( \, \lim_{t \to +\infty} r_t=0 \big )=\P \big (  \,\sup_{t \in \R_+} r_t<+\infty \big )=1\,.
\end{align*}
\item (Recurrence) If $\lambda_1 \ge 0 $ then $s(+\infty)=+\infty$,  $s(+0)=-\infty$. %and \\$\int_{(0,\infty)}m(x)\dd x<\infty $.
In this case we have
\begin{align*}
\P \big ( \, \inf_{t\in \R_+} r_t=0 \big )=\P \big ( \, \sup_{t \in \R_+} r_t=+\infty \big )=1\,.
\end{align*}
%and $r_t$ is positive recurrent, having invariant \textit{ probability} measure with density
%\begin{align*}
%\frac{1}{\int_0^\infty m(x)\dd x}m(x)\,.
%\end{align*}
%\item (Null Recurrence) If $\lambda_1 = 0$ then $s(+\infty)=+\infty$,  $s(+0)=-\infty$ and \\$\int_{(0,\infty)}m(x)dx=\infty $. In this case
%\begin{align*}
%\P \big (  \inf_{t< \infty} r_t=0 \big )=\P \big (  \sup_{t <\infty} r_t=+\infty \big )=1
%\end{align*}
%and $r_t$ is null recurrent, that is the process has no invariant \textit{probability} measure, and the mean time it takes the process, started at $x>0$ to come back is infinite.
\end{enumerate}
Since the speed measure is integrable around zero exactly when $\lambda_1>0$ and is always integrable around $+\infty$ we have that
$$\int_0^{+\infty} m(x)\,\dd x<\infty \,\,\, \Longleftrightarrow \,\,\,  \lambda_1>0$$
and therefore the measure
$$m_p(x)\dd x:=\frac{1}{\int_0^{+\infty} m(u)\,\dd u}m(x)\,\dd x$$
is a well defined probability measure which is invariant as $\mathcal A^\star m_p(x)=0$. In case $\lambda_1\le 0$ the speed measure is not finite and we do not have any invariant probability measures.
 \end{proof}
We finish this section by prooving the following lemma concerning the a.s. spatial regularity of IOUFs. The one-point-motion $\phi_t(x)$ of an IOUF is as already stated an Orntsein-Uhlenbeck-process and so if $x$ is far away from the origign $|\phi_t(x)|$ will decrease roughly as $|x|e^{-ct}$ (which  is the expected decrease). because the variance of $|\phi_t(x)|$ is rather negligible for small $t$. We may expect the IOUFs unitstep-discretisation $\phi=\phi_{0,1}$ to look like $e^{-c}$ time the Identity on a large scale. For a fixed $x$ this is quite obvious but care has to be taken about the fact that we are dealing with infinitely many random variables. The next lemma states that this is no problem at all.
\begin{lemma}\label{spat_reg_lem}
Let $\phi=\phi_{0,1}:\R^d\rightarrow\R^d$ be as in Definition~\ref{def_iouf} (with $s=0$ and $t=1$). Then we have a.s.
\begin{enumerate}
\item
\begin{equation}\label{eq_spacial}
\lim_{R\to\infty}\sup_{||x||\geq R}\frac{||\phi(x)-e^{-c}x|| }{||x||}=0
\end{equation}
\item \begin{equation}
       \lim_{R\to\infty}\sup_{||x||\geq R}\frac{||\phi(x)||}{||x||}=e^{-c}
      \end{equation}
\end{enumerate}
\end{lemma}
\begin{proof}
Observe that the second formula is an easy consequence of the first one, so we will only have to prove~\eqref{eq_spacial}. Therefore it is sufficient to show
\begin{equation}
 \lim_{R\to\infty}\sup_{R\le ||x||\le  R+1}\frac{||\phi(x)-e^{-c}x|| }{||x||}=0
\end{equation}
Let $A_R:=\{x\in\R^d:R\le ||x||\le  R+1 \}$. We want to apply the chaining technique (see Proposition~\ref{prop_chaining}) and so we first observe that $A_R$ can be covered by $c_3 R^d 3^j$ balls of radius $\delta_j=3^{-j}$. Let $\chi_j$ consist of the centers of these balls and $\epsilon_j=2^{-j-2}$ as well as $\epsilon=\frac{1}{2}$.  Fixing $\tilde{\epsilon}>0$ and an arbitrary $x_0\in A_R$ and assuming $R\ge \frac{16 e^{c}}{\tilde\epsilon}$ we conclude
\be
\prob{\sup_{x\in\A_R}\frac{||\phi(x)-e^{-c}x|| }{||x||}>\tilde\epsilon} \le\prob{\sup_{x\in\A_R}||\phi(x)-e^{-c}x|| >\tilde\epsilon R}
\le\prob{||\phi(x_0)-e^{-c}x_0|| >\frac{\tilde\epsilon R}{2}}\nonumber\\+\sum_{j=0}^\infty c_3R^d3^j\sup_{|x-y|\leq3^{-j}}\prob{|||\phi(x)-e^{-c}x||-||\phi(y)-e^{-c}y|||>2^{-j-2}\tilde\epsilon R}
\ee
Using the standard estimates for the normal distribution (as stated in Lemma~\ref{comp_lemma_2} for a Brownian motion) one easily shows that.
\[\prob{||\phi(x_0)-e^{-c}x_0|| >\frac{\tilde\epsilon R}{2}}\le c_4 e^{-\frac{\tilde\epsilon^2}{8d^2}R^2}\]
where $c_4$ is a constant only depending on $c$ and $d$. Using Corollary~\ref{iouf_gen_cor_0} and Lemma~\ref{comp_lemma_2} we get for $|x-y|\leq3^{-j}$:
\be
&&\prob{|||\phi(x)-e^{-c}x||-||\phi(y)-e^{-c}y|||>2^{-j-2}\tilde\epsilon R}\nonumber\\
&\le& \prob{|||\phi(x)-e^{-c}x||-||\phi(y)-e^{-c}y|||>2^{-j-2}\tilde\epsilon R3^j|x-y| }\nonumber\\
&\le&\prob{|||\phi(x)-\phi(y)|||>2^{-j-3}\tilde\epsilon R3^j|x-y| } \\
&\le& \prob{B_1^*\ge\frac{\log(2^{-3-j}\tilde\epsilon R3^j)-\lambda}{\sigma}}\le c_5 (2^{-j-3}\tilde{\epsilon R 3^j})^{-\frac{\log(2^{-3-j}\tilde\epsilon R3^j)-2\lambda}{2\sigma^2}}
\ee
where the constant $c_5$ depends only on $\sigma$ and $\lambda$. Combining all the above estimates we conclude that for arbitrary $\tilde\epsilon>0$ we have that $\prob{\sup_{x\in\A_R}\frac{||\phi(x)-e^{-c}x|| }{||x||}>\tilde\epsilon}$ is sumable over $R$ which yields the desired conclusion via an application of the Borel-Cantelli Lemma.
\end{proof}

\section{The  dimension of the statistical equilibrium of an isotropic Ornstein-Uhlenbeck flow }
\label{sec:dimsec}

Throughout this section  $\varphi\colon \R \times \R^d\times \Omega  \to \R^d $ will denote  an RDS over the  MDS $\left(\Omega,\F,(\theta_t)_{t\in \R},\P \right)$, associated, in the sense of Remark \ref{sf_rds}, to an isotropic Ornstein-Uhlenbeck stochastic flow $\phi\,\colon\, \R_+\times \R_+ \times \tilde \Omega \times \R^d \to \R^d$. We will also call $\varphi$ an IOUF.   The Lyapunov spectrum of $\varphi$ will be denoted by $\lambda_1>\dots>\lambda_d$. It has been calculated in Proposition \ref{loc_prop2}.\\[0.2cm]
\bem We also have the standing assumption that the top Lyapunov exponent is strictly positive.\eem\\[0.2cm]
Our goal here will be to calculate the dimension of the statistical equilibrium of $\varphi$.\\
In the next theorem we  quote some general facts  concerning the statistical equilibrium of an IOUF. In a slightly different setting it can be found in   \cite{LeJ85}, \cite{Bax91} and \cite{Ku90}, more precisely in \cite{LeJ85} for the case $\varphi$ is an IBF,  in \cite{Bax91} in the case $\varphi$ is generated by an SDE driven by  finitely many  Brownian motions. Virtually the same is also the proof of Theorem 4.3.6 in \cite{Ku90}, however the measures considered there  are given by $m_t(A):=m(\varphi(t,A,\omega))$ for positive $t$.
\begin{theo}\label{iouf_se_1}
Let $\varphi\,\colon\, \R \times \R^d\times \Omega  \to \R^d $ be an RDS over the  MDS $\left(\Omega,\F,(\theta_t)_{t\in \R},\P \right)$, associated to an isotropic Ornstein-Uhlenbeck stochastic flow $\phi$ %$\,\colon\, \R_+\times \R_+ \times \Omega \times \R^d \to \R^d$
 in the sense of Remark \ref{sf_rds}. Then the following statements hold:
\begin{enumerate}
\item[\textnormal{(a)}] There is an  invariant Markovian measure (see Proposition 1.4.3 in  \cite{Arn98}) $\mu$ on $(\Omega \times \R^d, \F\otimes \B(\R^d))$ with  factorization  $$\mu(\dd x,\dd\omega)=\mu_\omega(\dd x)\P(\dd\omega)$$
satisfying
$$\E^\P\mu_\omega=\nu\,,\,\,\, \text{ i.e. } \,\,\,\,\,\,\,\E\int f \dd \mu_\omega=\int f\dd \nu \text{ for all } f \in C_b(\R^d)\,,$$
where $\nu$ is the unique invariant probability measure for the one-point motion and $ C_b(\R^d) $ is the set of all bounded continuous functions mapping from $ \R^d $ to $ \R $.
\item[\textnormal{(b)}] The measure $\mu_\omega$ can be obtained as the $\P$-almost sure weak limit
$$\lim_{t \to \infty}\nu\circ \varphi^{-1}(t,\theta_{-t}\omega)=\mu_\omega\,.$$
\item[\textnormal{(c)}] For arbitrary $t\in \R$ it holds that
$$\mu_\omega \circ \varphi^{-1}(t,\omega)=\mu_{\theta_t \omega}\,\,\,\,\,\,\P\text{-almost surely}.$$
\end{enumerate}
\end{theo}
\begin{remark}
It is well known that
$$\nu(dx)=\left(\frac{c}{\pi}\right)^\frac{d}{2} e^{-c|x|^2}\,.$$
\end{remark}
Before we proceed with proving the main result in this section  we need to introduce the notion of a \begin{em} pointwise dimension \end{em} of a measure, which is in a sense a  refinement of the notion of density.
\begin{defi}
Let $\mu$ be a Borel measure on $\R^d$. The \bem lower and upper pointwise dimensions \eem are defined by
\begin{align*}
\underline{\dd}_\mu(x)=\liminf_{\epsilon \to+ 0}\frac{\log \mu(\mathrm{B}(x,\epsilon))}{\log \epsilon} \,\,\,\,\text{ and }\,\,\,\, \overline{\dd}_\mu(x)=\limsup_{\epsilon \to +0}\frac{\log \mu(\mathrm{B}(x,\epsilon))}{\log \epsilon}\,.
\end{align*}
If $\underline{\dd}_\mu(x)=\overline{\dd}_\mu(x)$ we call the common value the \bem  pointwise dimension \eem and denote it by $\dd_\mu(x)$. If there is a constant $\dd(\mu)=\dd_\mu(x)$ for $\mu$-almost all $x$ then $\dd(\mu)$ is called the dimension of the measure $\mu$.
\end{defi}
One can easily construct quite regular measures which do not have fixed dimension. For example  in $\R^2$ take a disjoint  ball $B$ and a line segment $L$ and place uniform distributions $U_B$ and $U_L$ on them. $U_B$ and $U_L$ have dimensions $2$ and $1$ respectively, but the dimension of $U:=U_B+U_L$ is not well defined.  Measures   having a fixed dimension are quite special and in a sense look similar around almost every point in a set supporting the measure. The main result in this chapter is that the statistical equilibrium of an IOUF is such a special measure having constant pointwise dimension. Further we give an explicit expression for the dimension using the results    of F. Ledrappier and L.-S. Young in \cite{LeY88}. There the authors    calculate the  dimension of the statistical equilibrium of a composition of independent and identically distributed diffeomorphisms on a compact space. There is no obstacle in the fact that we have an RDS in continuous time since the mappings $\left(\varphi(1,\theta_z\omega)\right)_{z \in \Z} $ are independent and identically distributed and define the same statistical equilibrium.\\
The only fact we still need to prove is that their assumption that the state space is compact is also not an obstacle since we have an invariant probability measure for the one-point motion and the RDS has a weak  attractor. This will be done in the course of the proof of the next theorem.   \\
One can link the dimension of a probability  measure $\mu$ on $\B(\R^d)$ (if it exists) to the minimal Hausdorff dimension of a set supporting  $\mu$:
\begin{align}\label{iouf_se_hdim}
\dd(\mu)=\inf\big\{\text{dim}_H Y : Y\in \B(\R^d),\, \mu(Y)=1\big\}\,,
\end{align}
where $ \text{dim}_H(Y) $ denotes the Hausdorff dimension of the set $ Y $.
For the proof of this fact see e.g.~\cite{LSY82}.
Observe that the above equality does not concern  \bem the support \eem of $\mu$ since we are not allowed to take the closure of the set $Y$. For example consider a probability measure $\mu$ on $\R$ with  $\mu([0,1]\cap \Q)=1$ supported on  $[0,1]$. The  Hausdorff dimension of the support is $1$, while in fact $\mu$ is carried by a countable set and thus has dimension zero.\\
We follow \cite{LeY88} and introduce the \bem Lyapunov dimension\eem:
\begin{defi}
Let $\lambda_1>\dots>\lambda_r$ be the Lyapunov exponents and $m_1,\dots,m_r$ the corresponding multiplicities of some RDS on a $d$-dimensional smooth ($C^\infty$) manifold. Let $k$ be the largest integer, such that $\sum\limits_{i=1}^k\lambda_i m_i>0$. The corresponding \bem Lyapunov dimension \eem $\D(\lambda_1,\dots,\lambda_r)$ is defined as follows:
\begin{align*}
\D(\lambda_1,\dots,\lambda_r)=\left \{  \begin{array}{ll}
 d & \text{if } \sum\limits_{i=1}^k m_i=d \\
  \sum\limits_{i=1}^k m_i-\frac{1}{\lambda_{k+1}}\sum\limits_{i=1}^k \lambda_i m_i& \text{otherwise}\,.
\end{array} \right .
\end{align*}
\end{defi}
\begin{theo}\label{iouf_se_mt}
Let $\varphi\,\colon\, \R \times \R^d\times \Omega  \to \R^d $ be an RDS over the MDS $(\Omega,\F, (\theta_t)_{t\in \R},\P)$ associated to an isotropic Ornstein-Uhlenbeck stochastic flow with strictly positive top Lyapunov exponent.  The  dimension of $\mu_\omega $ exists for almost all $\omega$, is deterministic  and one has
$$\mathrm{d}(\mu_\omega)=\mathcal{D}(\lambda_1, \dots,\lambda_d)\,.$$
\end{theo}
\textbf{Proof}:\\
Let us first  summarize the simple idea: Consider a one-point compactification of $\R^d$. We will associate a cocycle $\psi$ on the compactified space, which has the same Lyapunov spectrum as $ \varphi $. Further, the  dimension of the nontrivial statistical equilibrium of the new RDS will coincide with the one of $\varphi$. The  dimension of the statistical equilibrium of $\psi$ has been calculated in \cite{LeY88}. \\[0.3cm]
Let $\cal S$ be a $d$-dimensional sphere in $\R^{d+1}$ equipped with the induced topology from $\R^{d+1}$, and let $N$ be its north pole. Let $\mathcal S^N:=\mathcal{S} \setminus \{ N \}$.\\
It is well known that there is a smooth ($ C^\infty $) diffeomorphic map
$$g\colon \R^d\cup\{\infty\} \to \mathcal{S}\,\text{ such that } g(\R^d)=\mathcal{S} \setminus \{ N \}\,.$$
As usual we will denote by $T_pM$ the tangent space of the manifold $M$ at the point $p\in M$, by $TM$ the tangent bundle of the manifold $M$ and by $\Dd f$ the differential of a differentiable mapping $f\colon M \to N $, i.e $\Dd f\colon TM\to T N$ is given by
\begin{align*}
\Dd f(x,\cdot)\colon T_xM \to T_{f(x)}N,\hspace{0.4cm} (\Dd f)(x,v)=(f(x),\Dd f(x)v)\,.
\end{align*}
Consider the RDS $\psi^N$ on $\mathcal S^N$ over the metric DS $\left( \theta\right)_t$ given by
\begin{align*}
\psi^N(t,y,\omega)\colon \R  \times \mathcal S^N \times \Omega\to \mathcal S^N \hspace{0.4cm} y \to g\circ \varphi(t,\omega)\circ g^{-1}\,.
\end{align*}
That is, $ \psi^N $ is cohomologous to $ \varphi $ with cohomology $ g $ (see Proposition 1.9.6 on page 46 in \cite{Arn98}).
We will show that  $ \psi^N $ has the same Lyapunov spectrum as $\varphi$ by showing that their linearizations are Lyapunov cohomologous with cohomology $ Dg $ (see Proposition 4.1.9 in \cite{Arn98}).\\
According to Proposition 4.2.5 in \cite{Arn98}, the  linearization $T_1$ of $ \psi^N $ is a bundle RDS on $\Omega \times T\mathcal S^N$
\begin{align*}
T_1 \colon \Omega \times T\mathcal S^N\to \Omega \times T\mathcal S^N  \hspace{0.4cm} (\omega,(y,v)) \to (\theta_t \omega, \Dd \psi^N(t,\omega)(y,v))
\end{align*}
over the skew-product shift
$$\Theta_t\colon \Omega\times \mathcal S^N \to \Omega\times \mathcal S^N ,\hspace{0.5cm} \Theta_t(\omega,x)=(\theta_t \omega, \psi^N(t,\omega,y))\,.$$
$T_1$ is isomorphic (see Definition 1.9.8 in \cite{Arn98}) to the linearization $T_2$  of $\varphi$
\begin{align*}
T_2 \colon \Omega \times T \R^d\to \Omega \times T \R^d  \hspace{0.4cm} (\omega,(x,v)) \to (\theta_t \omega, \Dd \varphi(t,\omega)(x,v))
\end{align*}
via the mapping
\begin{align*}
K\colon \Omega \times T \R^d\to \Omega \times T\mathcal S^N\hspace{0.4cm}  (\omega,(x,v))\to (\omega, (g(x),\Dd g(x)v))
\end{align*}
since obviously $K(\omega,x,\cdot)\colon T_x \R^d \to T_{g(x)} \mathcal S^N$ is linear  on fibers. $T_2$ is a bundle RDS over the skew-product shift
$$\tilde\Theta_t\colon \Omega\times \R^d \to \Omega\times \R^d  ,\hspace{0.5cm} \Theta_t(\omega,x)=(\theta_t \omega, \varphi(t,x,\omega))\,.$$
Since $g\colon \R^d \to \mathcal S^N$ is diffeomorphic, the functions
\begin{align*}
&n_1 \colon \Omega \times \R^d \to (0,+\infty)\,, \hspace{0.4cm} (\omega,x)\to \log \|K(\omega,x)\|=\log\|\Dd g(x)\|\\[0.3cm]
&n_2 \colon \Omega \times \mathcal S^N \to (0,+\infty)\,, \hspace{0.4cm} (\omega,y)\to \log\|K^{-1}(\omega,y)\|=\log\|\Dd g^{-1}(y)\|\,,
\end{align*}
where $\|\cdot \|$ denotes the operator norm, are continuous. %Further we know that for all $x\in \R^d$
%\begin{align*}
%\varphi(t,x)\underset{t \to \pm\infty} \longrightarrow \mathcal{N}(0,2c\cdot\mathrm{Id}) \text{ in distribution}
%\end{align*}
%since $\varphi(t,x)|_{t \in \R_+}$ and $\varphi(-t,x)|_{t \in \R_+}$ are both Ornstein-Uhlenbeck processes.\\
%Therefore of course also for all  $y \in \mathcal S^N$
%\begin{align*}
%\psi(t,y)=g(\varphi(t,g^{-1}(y)))\underset{t \to \pm\infty} \longrightarrow \mathcal{N}(0,2c\cdot\mathrm{Id})\circ g^{-1} \text{ in distribution.}
%\end{align*}
%Combining this with the continuity of $n_1$ and $n_2$ we can follow
According to  Proposition 4.1.9 in \cite{Arn98}, in  order to show that $\psi^N$ has the same Lyapunov spectrum as $\varphi$ we need to show that for $\mu:=\mu_\omega(\dd x)\,\P(\dd \omega)$-almost all $(x,\omega)$ and respectively $\mu_\omega\circ g^{-1}(\dd y)\,\P(\dd \omega)$-almost all $(y,\omega)$
\begin{align}\label{iouf_se_mt1}
\begin{split}
&\lim_{t \to \pm\infty}\frac{1}{t}n_1(\tilde \Theta_t(\omega,x))=\lim_{t \to \pm\infty}\frac{1}{t}\log\|\Dd g(\varphi(t,x,\omega))\|=0\hspace{0.2cm} \,\,\, \text{  and}\\[0.4cm]
&\lim_{t \to \pm\infty}\frac{1}{t}n_2(\Theta_t(\omega,y))=\lim_{t \to \pm\infty}\frac{1}{t}\log \|\Dd g^{-1}(\psi^N(t,y,\omega))\|=0 \hspace{0.2cm} ,
\end{split}
\end{align}
which means that  $K$ is a \bem Lyapunov Cohomology \eem between the linearizations $T_1$ and $T_2$ (see Definition 4.1.6 in \cite{Arn98}).  According to Proposition 4.1.3 on page 165 of \cite{Arn98} (Dichotomy for linear growth of stationary processes) it is enough to show that $\mu_\omega(\dd x)\,\P(\dd \omega)$ and respectively $\mu_\omega\circ g^{-1}(\dd y)\,\P(\dd \omega)$-almost surely
\begin{align}\label{iouf_se_mt2}
\begin{split}
&\lim_{t \to +\infty}\frac{1}{t}n_1(\tilde \Theta_t(\omega,x))=\lim_{t \to +\infty}\frac{1}{t}\|\Dd g(\varphi(t,x,\omega))\|=0\hspace{0.2cm} \,\,\, \text{  and}\\[0.4cm]
&\lim_{t \to +\infty}\frac{1}{t}n_2(\Theta_t(\omega,y))=\lim_{t \to +\infty}\frac{1}{t}\|\Dd g^{-1}(\psi^N(t,y,\omega))\|=0.
\end{split}
\end{align}
However (\ref{iouf_se_mt2}) holds true and the main reason is that the invariant measure $\mu$ has Markov nature, i.e.~the factorization $\mu_\omega $ is measurable with respect to the past ($\F_{-\infty}^0$-measurable).
First observe that for every fixed $x\in \R^d$ and $y \in \mathcal{S}^N$
\begin{align*}
\varphi(t,x)\underset{t \to \pm\infty} \longrightarrow \mathcal{N}(0,\frac{1}{2c}\cdot\mathrm{Id}) \text{ in distribution}
\end{align*}
which follows from the convergence of an Ornstein-Uhlenbeck process to its equilibrium and
\begin{align*}
\psi^N(t,y)=g(\varphi(t,g^{-1}(y)))\underset{t \to \pm\infty} \longrightarrow \mathcal{N}(0,\frac{1}{2c}\cdot\mathrm{Id})\circ g^{-1} \text{ in distribution}
\end{align*}
because convergence in distribution is preserved by continuous mappings. Combining these with the continuity of $\log \|\Dd g(\cdot)\|$ and $\log \|\Dd g^{-1}(\cdot)\|$ implies that the convergences (\ref{iouf_se_mt2})  hold for every fixed $x\in \R^d$ and $y \in \mathcal{S}^N$, $\P$-almost surely. Now, since $\mu_\omega$ is independent of $\F_0^t$ for every $t>0$ we obtain that (\ref{iouf_se_mt2}) holds also for $\mu_\omega(\dd x)\,\P(\dd\omega)$-almost all $ (x,\omega) $ and respectively for $\mu_\omega\circ g^{-1}(\dd y)\,\P(\dd\omega)$-almost all $ (y,\omega) $. Therefore $ K $ is a Lyapunov cohomology and the Lyapunov spectrum of $\psi^N$ is the same as the one of $\varphi$.\\
Now define the RDS $\psi$
\begin{align*}
\psi(t,y,\omega)=\left \{  \begin{array}{ll}
\psi^N(t,y,\omega) & \text{if }  y \ne N, \\
 N & \text{if }  y = N.
\end{array} \right .
\end{align*}
It is an  RDS on $\mathcal S$ and its linearization is a bundle RDS over the skew-product metric DS $\Theta_t(\omega, y)=(\theta_t \omega, \psi(t,y,\omega))$, preserving the measure $\mu_\omega\circ g^{-1}(\dd y)\,\P(\dd\omega)$. However, since this measure puts zero mass on the $\Theta_t$-invariant set $\{N\}\times \Omega$, and $\psi^N$ differs from $\psi$ only on this set, we can conclude that the Lyapunov spectrum of the linearization of $\psi^N$ coincides with the Lyapunov spectrum of the linearization of $\psi$, viewed as a bundle RDS over the metric DS $\Theta_t$ preserving $\mu_\omega\circ g^{-1}(\dd y)\P(\dd\omega)$.\\
Despite the fact that one can use the Lemma~\ref{spat_reg_lem} to
show the differentiability of $\psi$, it is easier to consider the two
ergodic components of $\psi$, each of them being a $C^2$ random dynamical system. Once started in $\mathcal{S}\setminus N$ one stays there and hence the differentiability in $N$ is not important for our purposes, since the flow couldn't even
observe a lack of it.
%Note that the linearization of $\psi$ can also be viewed as a bundle RDS over  the skew-product metric DS $\Theta_t(\omega, y)=(\theta_t \omega, \psi(t,\omega,y))$, preserving the measure $\delta_N(\dd y)\P(\dd\omega)$ and in this case all Lyapunov exponents are $0$. However we are not interested in this point of view.\\

Now we want to apply the results of Ledrappier and Young \cite[Theorem A]{LeY88} giving the dimension of the statistical equilibrium of an  ergodic  random dynamical system  and therefore to be in their setting we will consider the discretized cocycle $\psi\colon \Omega\times\R^d\times \Z \to \R^d$. In our case the assumptions in \cite{LeY88} are  verified: \\
1. Clearly the shift $\theta:=\theta_1$ is ergodic.  \\
2. The set of invariant probability measures of the one-point motion of $\psi$ is given by $\{\lambda \delta_N+(1-\lambda)\mathcal{N}(0,\frac{1}{2c}\cdot\mathrm{Id})\circ g^{-1} : \lambda \in [0,1] \}$ and thus $\nu\circ g^{-1}:=\mathcal{N}(0,\frac{1}{2c}\cdot\mathrm{Id})\circ g^{-1}$ is an extremal point in this convex set.\\
Further, the technical assumptions
\begin{align*}
\E \log^+||\psi(1)||< \infty \text{ and } \E \log^+||\psi(1)^{-1}||< \infty
\end{align*}
are direct consequences of Theorem 2.1(v) from \cite{MoSch99}, after taking into account that $\psi$ acts on the compact space $\mathcal{S}$.\\
The last condition which  has to be checked is that one of the Hypothesis $A$, $A'$ or $B$ from \cite{LeY88} is satisfied. It can be done for example using the remarks on stochastic flows in Section 5.2  of \cite{LeY88} by showing that the diffusion $(\varphi(t,x), \Dd \varphi(t,x)e)$ is hypoelliptic  and therefore the same holds for  $\psi$.\\
We will check however the \bem Hypothesis \eem $A$:\\
3. The MET gives  splitting of $T_y\mathcal S$ as direct sum $E_1(y,\omega)\oplus\dots \oplus E_m(y,\omega)$ of linear subspaces corresponding to the distinct Lyapunov exponents $\lambda_1 >\dots>\lambda_d$. The hypothesis states that for $\nu\circ g^{-1}$-almost all $y$ the distribution of
$\omega \to E^j(y,\omega)$ for $j=k$ and $j=k+1$, where $k$ is as in the definition of the Lyapunov dimension, is absolutely continuous with respect to the Riemannian measure on the space of $\sum_{i\ge j}m_i$-planes in $T_y\mathcal{S}$. Here $m_i$ is the multiplicity of $\lambda_i$ and  $E^j(y,\omega):=\oplus_{i\ge j }E_i(y,\omega)$.\\
If we can verify the hypothesis for $\varphi$ it will be also true for $\psi^N$ on $\mathcal{S}^N$, because of Proposition 4.1.9 in \cite{Arn98} and the fact that   $\varphi$ and $\psi^N$ are Lyapunov cohomologous. Further it will also hold for $\psi$ as $\nu\circ g^{-1}\otimes\P(\delta_N\times \Omega)=0$.
However, for arbitrary orthogonal transformation $O\in \mathcal O_d$ the isotropy   of $\varphi$ implies that the distribution of $ O\varphi(t,\cdot)$ coincides with the distribution of $\varphi(t, O \cdot)$ and therefore the distributions of $ O \Dd \varphi(t,0)$ is the same as the one of $\Dd\varphi(t,0)$, and thus for the splitting $E_1(0,\omega)\oplus\dots \oplus E_d(0,\omega)$ of $T_0\R^d$ we have that the distribution of $\omega \to E^j(0,\omega)$ is exactly uniform on the space of $(d-j)$-planes in $T_0 \R^d$. Actually, this holds also for arbitrary $x\in \R^d$ since the distribution of the semimartingale field
$M_x(t,y):=\int_0^t \Dd F(s,\varphi(s,x))y-cy $ generating  $ \Dd \varphi(t,x) $ via an SDE is independent of $x$ and therefore coincides with the distribution of the semimartingale field  generating $\Dd\varphi(t,0)$ (see (\ref{iouf_gen_prop2.1}) and the discussion before  in the proof of Proposition \ref{loc_prop2}).\\
The preceding argument  depends heavily on the rotation invariance of the generating field. Much more general treatment on the laws of the Osceledets spaces of cocycles arising as solutions of linear SDE is provided by Peter Imkeller in \cite{Imk97} and  \cite{Imk98}. \\
Now we use~\cite[Theorem A]{LeY88} and state that the  dimension of $\mu_\omega\circ g^{-1}$ exists $\P$-almost surely and
\begin{align*}
\text{d} (\mu_\omega\circ g^{-1})=\mathcal{D}(\lambda_1,\dots,\lambda_d)\,.
\end{align*}
It remains to be shown  that the  dimension of $\mu_\omega$ also exists and  equals the  dimension of $\mu_\omega\circ g^{-1}$. This is done via the following:\\
\textit{Claim}: Let $(X,\dd_X)$ and $(Y,\dd_Y)$ be metric spaces and let $f\colon X \to Y$ be a locally Lipschitz bijection, such that $f^{-1}$ is also locally Lipschitz. Let further $\mu$ be a compactly supported positive measure on $(X,\B(X))$ with $\dd(\mu)=d\,.$
Then the  dimension of the image measure $\mu\circ f^{-1}$ exists and equals $d$.\\[0.3cm]
\textit{Proof of the Claim}:\\
Let $K >0$ be such that
$$\dd_Y(f(x_1),f(x_2))\le K \dd_X(x_1,x_2) \text{ and }\dd_X(f^{-1}(y_1),f^{-1}(y_2))\le K \dd_Y(y_1,y_2)$$
for all $x_1,x_2 \in \text{supp}(\mu)$ and $y_1,y_2\in \text{supp}(\mu\circ f^{-1})$.
Let $\text{B}_Y(\tilde y,R)$ be a ball in $Y$ with  radius $R$ and centered at $\tilde y\in Y$.
Then obviously
\begin{align*}
f^{-1}(\text{B}_Y(\tilde y,R))\subset \text{B}_X(f^{-1}(\tilde y),K R)
\end{align*}
This implies
\begin{align*}
&\limsup_{R\to 0}\frac{\log \mu\circ f^{-1}(\text{B}_Y(\tilde y,R))}{ \log R}=\limsup_{R\to 0}\frac{\log \mu \left (f^{-1}(\text{B}_Y(\tilde y,R))\right )}{\log R}\le \limsup_{R\to 0}\frac{\log \mu \left (\text{B}_X(f^{-1}(\tilde y),KR)\right )}{\log R}\\
&= \limsup_{R\to 0}\frac{\log \mu \left (\text{B}_X(f^{-1}(\tilde y),KR)\right )}{\log (K R)-\log K}=d \eins_{\text{supp}(\mu)}(f^{-1}(\tilde y))=d \eins_{\text{supp}(\mu\circ f^{-1})}(\tilde y)\,.
\end{align*}
On the other hand
\begin{align*}
\text{B}_{X}(f^{-1}(\tilde y),\frac{1}{K}R)\subset f^{-1}(\text{B}_Y(\tilde y,R))
\end{align*}
and therefore
\begin{align*}
&\liminf_{R\to 0}\frac{\log \mu\circ f^{-1}(\text{B}_Y(\tilde y,R))}{ \log R}=\liminf_{R\to 0}\frac{\log \mu \left (f^{-1}(\text{B}_Y(\tilde y,R))\right )}{\log R}\ge \limsup_{R\to 0}\frac{\log \mu \left (\text{B}_X(f^{-1}(\tilde y),\frac{1}{K}R)\right )}{\log R}\\
&= \liminf_{R\to 0}\frac{\log \mu \left (\text{B}_X(f^{-1}(\tilde y),\frac{1}{K}R)\right )}{\log (\frac{1}{K} R)+\log K}=d \eins_{\text{supp}(\mu)}(f^{-1}(\tilde y))=d \eins_{\text{supp}(\mu\circ f^{-1})}(\tilde y)\,.
\end{align*}
Combining these two bounds for the $\limsup$ and the $\liminf$ we obtain
\begin{align*}
\lim_{R\to 0}\frac{\log \mu\circ f^{-1}(\text{B}_Y(\tilde y,R))}{ \log R}= d \eins_{\text{supp}(\mu\circ f^{-1})}(\tilde y)\,,
\end{align*}
that is the local pointwise dimension of $\mu\circ f^{-1}$ exists and equals the one of $\mu$. \\[0.3cm]
Now to complete the proof of the theorem we just have to notice that the support of the measure $\mu_\omega$ is almost surely compact, since the support is contained in the global weak set  attractor (see \cite{Gu449}, Theorem 2),  which is  compact. The existence of a global weak set attractor for an IOUF is obtained in the next section. Of course this implies that $\mu_\omega\circ g^{-1}$ has a compact support in $\mathcal{S}^N$. Now apply the statement of the claim $\omega$-wise to obtain that $\P$-almost surely
\begin{align*}
\dd (\mu_\omega) =\mathcal{D}(\lambda_1,\dots,\lambda_d)  \,,
\end{align*}
which completes the proof.\hfill
$\square$\\[0.5cm]
\textbf{Theorem~\ref{iouf_se_mt}} might be of interest also because one can link the support of the statistical equilibrium to the global set and point attractors, as it is done in \cite{Cr01} and \cite{KuSh04}. \\
Theorem 4.3 in \cite{Cr01} states that if $\varphi$ is a white noise RDS, that is it has independent increments, then the strong global point attractor $ A_p^s(\omega) $, if it exists,  supports every invariant probability measure , i.e. for the factorization $\mu_\omega$ of an arbitrary invariant measure it holds
$$\mu_\omega(A_p^s(\omega))=1.$$
Clearly we have also $\mu_\omega(A^s(\omega))=1$, where $A^s(\omega)$ is the global strong attractor, as in general $A_p^s(\omega)\subset A^s(\omega)$. \\
Gunter Ochs has shown in \cite{Gu449} (Theorem 2) that $ \mu_\omega(A_p(\omega))=1 $ where $ A_p(\omega) $ is the minimal global weak point attractor if it exists.
Kuksin and Shirikyan  \cite{KuSh04} (Proposition  1.6)   consider the relations between the weak global point attractor and the support of the random measure $\mu_\omega$. They consider RDS with discrete time  and independent increments. They  have shown that under the conditions C1 and C2 stated below (additional to the assumption of independent increments),  the support of the factorization $\mu_\omega$ obtained via  a limiting procedure as in Theorem~\ref{iouf_se_1} coincides almost surely ($\mu_\omega$ is itself $\P$-almost surely defined!) with the unique minimal weak point attractor if it exists. The conditions C1 and C2 are:
\begin{enumerate}
\item[C1] The one point motion has unique invariant probability measure $\nu$ and the family $(\varphi(t,x))_{t\in \R}$ converges in distribution as $t\to \infty$ to this invariant probability $\nu$ for all $x$.
\item[C2] For any $x\in \R^d$ and $\epsilon>0$, there is a $\Omega_\epsilon \in \F$, a compact set $K_\epsilon\subset \R^d$ and an integer $k_\epsilon$, such that $\P(\Omega_\epsilon)\ge 1-\epsilon$ and
$$\varphi(k,x,\theta_{-k}\omega)\in K_\epsilon \,\,\,\,\,\,\text{ for all } \omega \in \Omega_\epsilon \text{ and } k\ge k_\epsilon\,.$$
\end{enumerate}
Observe that  the condition C2 is not automatically satisfied if the RDS has a weak point attractor, since we need the set $\Omega_\epsilon$ to be independent of $k$. However, even if we know that the weak point attractor coincides with the support of the statistical equilibrium we still get only a bound from below for its Hausdorff dimension via the dimension of the statistical equilibrium, because the latter provides a bound from below for  the Hausdorff dimension of $\mathrm{supp}(\mu_\omega)$ (recall (\ref{iouf_se_hdim}) and the discussion below).

\begin{propo}
Let $\varphi$ be a RDS on $\R^d$ with time $\R$ over the MDS $(\Omega,\F,(\theta_t)_{t\in \R},\P)$ with global weak  point attractor $A_p(\omega)$. Let also $\varphi$ correspond to some IOUF $\phi$ with strictly positive Lyapunov exponent in the sense that $\{\phi_{s,t} : s,t\ge 0\}$ and $\{\varphi(t-s,\theta_s \omega) :  s,t\ge 0\}$ coincide in distribution. Then
\begin{align*}
\textnormal{dim}_H (A_p(\omega)) \ge \mathcal{D}(\lambda_1,\dots,\lambda_d)  \,,
\end{align*}
where $\textnormal{dim}_H(A)$ denotes the Hausdorff dimension of  $A$.
\end{propo}
\textbf{Proof:}\\
As mentioned above, the existence of a weak point attractor $A_p(\omega)$ is supplied by Theorem \ref{iouf_atr_main} and we take (as in \cite{KuSh04}) the minimal one. The fact that $\text{supp}(\mu_\omega)\subset A_p(\omega)$ almost surely has been discussed above. The rest is due to the fact that $\dd(\mu)=\inf\big\{\text{dim}_H Y : Y\in \B(\R^d),\,  \mu(Y)=1\big\}$ (see e.g. \cite{LSY82}). \hfill $\square$

We expect that the dimension of the statistical equilibrium of an IBF can be obtained as the limit of the dimensions of the statistical euilibrium of the corresponding IOUF as $c\to 0$. In the very special volume preserving case where the Lyapunov exponents of the IBF sum up to zero one easily sees that indeed $\mathcal{D}_c(\lambda_1,\dots,\lambda_d)\to d$.

\section{Weak attractors for isotropic Ornstein-Uhlenbeck flows}

\subsection{Weak attractors}
In this section we give a brief introduction to the notion of a weak random attractor.  Recall that
 the \bem Hausdorff distance \eem $\dd_H(A,B)$ between two  subsets $A,\,B\subset X$  of a metric space $(X,d)$ is defined as  $\dd_H(A,B)=\sup_{x \in A} \inf_{y \in B}d(x,y)$.
\begin{defi}\label{f_ibf_rcs}
A random  subset $A(\omega)\subset \R^d$ on the probability space $(\Omega,\F,\P)$ is called a \begin{em}random compact set\end{em} if
\begin{enumerate}
\item[\textnormal{(a)}] $A(\omega)$ is   nonempty and compact for all $\omega\in \Omega$.
\item[\textnormal{(b)}] The mapping $\omega \mapsto \mathrm{dist}(x,A(\omega))$ is measurable for each $x\in \R^d$.
\end{enumerate}
\end{defi}
We now give  a definition of a global weak attractor following \cite{Gu449}.
\begin{defi}\label{f_ibf_watr}
Let $\varphi$ be a RDS on $\R^d$ over  the metric DS $(\Omega,\F,\left(\theta_t\right)_{t\in \R},\P)$. The random set $A(\omega)$ is called a \bem weak global  set attractor \eem for $\varphi$ if
\begin{enumerate}
\item[\textnormal{(a)}] $A(\omega)$ is a random  compact set.
\item[\textnormal{(b)}] $A$ is strictly $\varphi$-invariant, that is $\varphi(t,\omega)(A(\omega))=A(\theta_t\omega)$ for  all $\omega \in \Omega$ and $t\in\R$.
\item[\textnormal{(c)}] $\lim\limits_{t\to \infty}\dd_H(\varphi(t,\theta_{-t}\omega)(B), A(\omega))=0$  in probability for all compact sets $B$.
\end{enumerate}
\end{defi}
%\begin{iouf_atr_rem}
%Since we will consider only RDS on $\R^d$ in Definition \ref{f_ibf_watr} (c) it is enough to require the almost sure convergence  for all  compact subsets.
%\end{iouf_atr_rem}
\begin{remark}
Actually the original definition in \cite{Gu449}  requires that the weak attractor attracts all random sets in some basin. However, in our special case (global set attractor, state space $\R^d$) the requirement of attracting all deterministic compact sets already determines the weak set attractor uniquely (see Theorem 3 and Corollaries 3.1 and 3.2 in \cite{Gu449}).
\end{remark}
Observe that the random attractor of a RDS gives comprehensive information on the asymptotics of the flow.\\
If in Definition  \ref{f_ibf_watr} (c) we replace the  convergence in probability with almost sure convergence we obtain the \bem global strong set attractor \eem as introduced  in \cite{CrFl94}.
\begin{remark}
According to \cite{Cr01}(Remark 3.2(ii)) or \cite{Cr99}(Section 5) the global strong set attractor is unique. It is already  uniquely determined by attracting all compact subsets.
\end{remark}
From now on, if not explicitly stated, we will always consider global weak set attractors. \\
 One can also define  the global weak and strong point attractors:
\begin{defi}\label{f_ibf_patr}
Let $\varphi$ be a RDS on $\R^d$ over  the metric DS $(\Omega,\F,\left(\theta\right)_{t\in \R},\P)$. The random set $A_p(\omega)$ is called a \bem global strong (weak) point   attractor \eem for $\varphi$ if
\begin{enumerate}
\item[\textnormal{(a)}] $A_p(\omega)$ is a random  compact set.
\item[\textnormal{(b)}] $A_p$ is strictly $\varphi$-invariant, that is $\varphi(t,\omega)(A_p(\omega))=A_p(\theta_t\omega)$ for  all $\omega \in \Omega$ and $t\in\R$.
\item[\textnormal{(c)}] $\lim\limits_{t\to \infty} \text{dist}(\varphi(t,x,\theta_{-t}\omega), A_p(\omega))=0$  $\P$-almost surely (in probability) for all $x\in \R^d$.
\end{enumerate}
\end{defi}
As it has been already mentioned in \cite{Cr01}  point attractors need not be unique. However it is a direct consequence of  Theorem 3.4 in \cite{Cr01} that if a global strong point attractor exists,  then there is a unique minimal point attractor.  Clearly the global strong set attractor is also a point attractor but not necessarily the minimal one. The same is stated for weak point attractors in \cite{KuS04}.\\[0.5cm]

\subsection{The case of IOUF}
Random attractors are defined for two-sided RDS, however our object of study is a class of stochastic flows. As already mentioned in Remark \ref{sf_rds} to a stochastic flow $\phi$ with stationary increments generated via a Kunita type SDE and satisfying some regularity conditions one can canonically construct a corresponding two-sided RDS $\varphi$ having the same distribution.
When we say that the IOUF has a weak attractor we mean that the corresponding RDS has this property.
We now state the main result in this section.
\begin{theo}\label{iouf_atr_main} Every isotropic Ornstein-Uhlenbeck flow  $\phi$  has a global weak set attractor.
\end{theo}

The proof of the theorem will rely on the following general criteria for weak attractors, which might be of interest also in a different context. An ongoing joint work of one of the authors with Michael Scheutzow and Hans Crauel  \cite{CDS08} deals among others with the following criteria in a very general framework. The proof in the following rather restricted setting can be found in \cite{Dim06}.
\begin{theo}
Let $\varphi\colon \R\times \R^d\times \Omega \to \R^d $  be a perfect cocycle on $\R^d$ over the metric dynamical system $\big( \Omega,\F,(\theta_t)_{t \in \R},\P\big)$ . Then the following are equivalent:
\begin{itemize}
\item[\textnormal{(i)}]There is a random DS,  $\tilde \varphi$ over the same MDS, which is indistinguishable from $ \varphi $\\ $($i.e.~$ \P(\varphi(t)=\tilde\varphi(t) \text{ for all } t\in \R )=1 )$ and  has a global weak set  attractor.
\item[\textnormal{(ii)}]For every  $\epsilon>0$ there exists $R_0>0$, such that for all $R>0$ there is a $t_0>0$  with the property that  for all $t\ge t_0$,
\begin{align*}
\P\big(\varphi(t)(\textnormal{B}(0,R))\subset \textnormal{B}(0,R_0)\big)\ge 1-\epsilon\,.
\end{align*}
\end{itemize}
\end{theo}

We will actually not verify the sufficient and necessary condition for the existence of a  weak random attractor stated above but rather use the following simple corollary of the preceding theorem. The proof is a straightforward application of the Markov inequality.
\begin{corollary}\label{iouf_atr_cor}
Let $\varphi$  be an RDS on $ \R^d $ over the MDS $ (\Omega,\F,(\theta_t)_{t\in \R},\P) $. Assume that there is a constant $M>0$, such that  for all $R>0$ there exists  $t_0:=t_0(R)>0$ with the property that  for all $t \ge t_0$
\begin{align*}
\E \big[\sup\{|\varphi(t,x)| \hspace{0.1cm} \colon  x\in \mathrm{B}(0,R)\}\big]\le M \textrm{ .}
\end{align*}
Then $\varphi$ has an indistinguishable  version $\tilde \varphi$ in the sense of the previous theorem, which has a global weak set attractor.
\end{corollary}

\subsection{Proof of the existence of weak attractor}
The proof of the main theorem depends crucially on the \bem chaining technique\eem..
We will use the formulation  from \cite{CSS00}, where also the proof can be found.\\
Let $(\X,d)$ be a compact metric space and $\phi \colon \X \to \R_+$ be a random continuous function i.e. a random variable taking values in the set of continuous functions from $\X$ to $\R_+$. Given a sequence of positive real numbers $\left(\delta_i\right)_{i \ge 0}$, such that $\sum\limits_{i=0}^\infty \delta_i <\infty$ we determine a sequence $\left(\chi_i\right)_{i=0 }^\infty$  of discretizations (skeletons) of $\X$, with the property that for all $x\in \X$ there is a point $x_i \in \chi_i$, such that $d(x,x_i)\le \delta_i$. Assume that $\chi_0=\{x_0\}$, with $d(x,x_0)\le \delta_0$ for all $x\in \X$.
\begin{propo}(Chaining)\label{prop_chaining}
Let $\phi\colon \X \to \R_+$ be an almost surely continuous random function with  $\left( \delta_i \right )_{i \ge 0}$ and $\left( \chi_i \right )_{i \ge 0}$  as above. For arbitrary positive $\epsilon,z\ge  0$ and an arbitrary sequence of positive reals $\left(\epsilon_i\right)_{i \ge 0}$ such that $\epsilon+\sum\limits_{i=0}^\infty \epsilon_i= 1$ we have
\begin{align*}
\P \big ( \sup_{x \in \X} \phi(x)>z\big)\le \P\big (  \phi(x_0)>\epsilon z\big)+\sum_{i=0}^\infty |\chi_{i+1}|\sup_{d(x,y)\le \delta_i}\P\big(  |\phi(x)-\phi(y)|>\epsilon_i z\big)\,.
\end{align*}
\end{propo}
\begin{proof}: See \cite[Lemma~4.1]{CSS00}.\end{proof}

The following lemma simply states some well known facts about the running maximum of a standard Brownian motion as well as a common estimate for the Gaussian tails, which we  use often below and therefore are stated explicitly.
\begin{lemma}\label{comp_lemma_2}
Let $\left (B_t \right )_{t \ge 0}$ be a standard Brownian motion, and let $B_t^\star:=\sup_{s\le t}B_s $ be its running maximum. The distribution of $B_t^\star$ has density $$\eins_{[ 0,\infty)}(x)\sqrt{\frac{2}{\pi t}}e^{-\frac{x^2}{2t}}$$ with respect to the Lebesgue measure.
Moreover for arbitrary $c>0$ the following bounds hold:
\begin{align*}
\P(B_t\ge c)\le \frac{1}{c}\sqrt{ \frac{t}{2 \pi }} e^{-\frac{c^2}{2t}} \hspace{1cm}\text{ and }\hspace{1cm} \P (B_t^\star\ge c)\le \frac{1}{c}\sqrt{ \frac{2 t}{ \pi }} e^{-\frac{c^2}{2t}}\,.
\end{align*}
\end{lemma}

The proof of the main result will be splited into several lemmas.

\begin{lemma}\label{O_U_estimate}
For all $\gamma_0>0$ there exist $t_0>0$ and $R_0>0$, such that for all $\gamma\ge \gamma_0$, $t \ge t_0$ and $R \ge R_0$
\begin{align*}
\P(|\phi_t(x)|>\gamma R\vspace{1cm})\le  \exp(-k \gamma^2 R^2)
\end{align*}
holds for all $x$ with $|x|=R$, where $k$ is a strictly positive constant depending on  the dimension of the state space $d$ and the slope of the potential $c$.
\end{lemma}

\begin{proof}
Because of the rotational invariance of the law of $\phi_t$ we can assume without loss of generality that $x=(R,0,\dots,0)^T$. Observe that
\begin{align*}
&\P\big(|\phi_t(x)|\ge \gamma R\big)=\P\big(\sum_{i=1}^d\phi_t^i(x)^2\ge \gamma^2 R^2\big)\\
&\le \P \big(\phi_t^1(x)^2\ge \frac{1}{4}\gamma^2 R^2\big)+\sum_{i=2}^d\P \big(\phi_t^i(x)^2\ge \frac{3}{4(d-1)}\gamma^2 R^2\big)\\
&=\P \big(|\phi_t^1(x)|\ge \frac{1}{2}\gamma R\big)+\sum_{i=2}^d\P \big(|\phi_t^i(x)|\ge \gamma R\sqrt \frac{ 3}{4(d-1)}\big)\,.
\end{align*}
Since $\phi_t(x)$ is an Ornstein-Uhlenbeck process started in $(R,0,\dots,0)$, the distributions of $\phi_t(x)^1$ and $\phi_t(x)^i$ for $i=2,\dots,d$ are  given by   $\mathcal N(R e^{-ct},\frac{1-e^{-2ct}}{2c})$ and  $\mathcal N(0,\frac{1-e^{-2ct}}{2c})$  respectively (see \cite{KaShr88}, section 5.6, page 358). Therefore, with  a standard Gaussian random variable $B$, we have:
\begin{align*}
&\P \big(|\phi_t(x)|\ge \gamma R\big)\\[0.4cm]
&\le \P\Big( |\sqrt \frac{1-e^{-2ct}}{2c}B+R e^{-ct}|\ge \frac{1}{2}\gamma R\Big)+\sum_{i=2}^d\P\Big(|\sqrt \frac{1-e^{-2ct}}{2c}B|\ge \gamma R\sqrt \frac{ 3}{4(d-1)}\Big)\\[0.4cm]
&=\P \Big( B \notin (R\sqrt \frac{2c}{1-e^{-2ct}} (-\frac{\gamma}{2}-e^{-ct})\,,\,R\sqrt \frac{2c}{1-e^{-2ct}} (\frac{\gamma}{2}-e^{-ct}) )\Big)\\[0.4cm]
&+\sum_{i=2}^d\P \Big(B \notin ( -\gamma R \sqrt \frac{ 3}{4(d-1)}\sqrt \frac{2c}{1-e^{-2ct}}\,,\, \gamma R \sqrt \frac{ 3}{4(d-1)}\sqrt \frac{2c}{1-e^{-2ct}}\Big)\\[0.4cm]
&\le 2 \P \Big ( B \ge R\sqrt \frac{2c}{1-e^{-2ct}} (\frac{\gamma}{2}-e^{-ct}) \Big )+2\sum_{i=2}^d\P \Big (B \ge  \gamma R \sqrt \frac{ 3}{4(d-1)}\sqrt \frac{2c}{1-e^{-2ct}}\Big )\\[0.4cm]
&\le 2 \P \big( B \ge R\sqrt{2c} (\frac{\gamma}{2}-e^{-ct}) \big )+2\sum_{i=2}^d\P \Big (B \ge  \gamma R \sqrt \frac{ 6c}{4(d-1)} \Big )\,.
\end{align*}
Choose $t_0>0$ large enough for $\frac{\gamma}{2}-e^{-ct}\ge \frac{\gamma}{4}$ to hold for all $t \ge t_0$. Thus for all $t \ge t_0$ with the help of \\Lemma \ref{comp_lemma_2} we obtain
\begin{align*}
&\hspace{0.6cm}\P\big(|\phi_t(x)|\ge \gamma R\big)\le 2 \P\Big( B \ge R\sqrt{2c} \frac{\gamma}{4} \Big )+2\sum_{i=2}^d\P\Big (B \ge  \gamma R \sqrt \frac{ 6c}{4(d-1)} \Big )\\[0.4cm]
&\le 2 \frac{1}{\sqrt{2\pi}}\frac{1}{R\sqrt{2c}\frac{\gamma}{4}}\exp \Big (  -\frac{1}{2}\frac{2c}{16}R^2 \gamma^2\Big )+2 \frac{d-1}{\sqrt{2\pi}}\frac{1}{\gamma R \sqrt\frac{6c}{4(d-1)}}\exp \Big (  -\frac{1}{2}\frac{6c}{4(d-1)}R^2 \gamma^2\Big )\,.
\end{align*}
Let $R_0>0$ be such that
\begin{align*}
\frac{2}{\sqrt{2\pi}}\frac{1}{R_0\sqrt{2c}\frac{\gamma_0}{4}}\le \frac{1}{2} \hspace{0.5cm}\textrm{ and } \hspace{0.5cm}2 \frac{d-1}{\sqrt{2\pi}}\frac{1}{\gamma_0 R_0 \sqrt\frac{6c}{4(d-1)}} \le \frac{1}{2} \,.
\end{align*}
Then for all $R\ge R_0$ and $\gamma \ge \gamma_0$ (and of course $t \ge t_0$) and with $k=\min\{\frac{c}{16},\frac{3c}{4(d-1)}\}$ we obtain
\begin{align*}
&\P\big(|\phi_t(x)|\ge \gamma R\big)
&\le \frac{1}{2}\exp \Big (  -\frac{c}{16}R^2 \gamma^2\Big )+\frac{1}{2}\exp \Big (  -\frac{1}{2}\frac{6c}{4(d-1)}R^2 \gamma^2\Big ) \le 2 \exp\big(-k \gamma^2 R^2\big)\,,\,\,
\end{align*}
which proves the statement.
%\hfill $\square$\\[0.5cm]
%The proof follows immediately from the fact that the one point motion $\phi_t(x)$ is an Ornstein-Uhlenbeck process and therefore subject to a Gaussian distribution. %The complete argument can be found in \cite{Dim06}
\end{proof}

The following lemma uses the chaining technique and is  crucial for the proof of the main theorem.
\begin{lemma}\label{iouf_atr_lem3}
For every $t>0$  there exist constants $K>0$ and $L>0$ (depending on $t$!), such that for all $R>0$
\begin{align*}
\E \big[\sup\{\sup_{s \in [0,t]}|\phi_{s}(x)| \hspace{0.1cm} \colon  x \in \textnormal{B}(0,R)\}\big]\le KR+L \,.
\end{align*}
\end{lemma}
\begin{proof}
Observe the simple computation:
\begin{align*}
&\E \big[\sup\{\sup_{s \in [0,t]} |\phi_{s}(x)| \hspace{0.1cm} \colon  x \in \mathrm{B}(0,R)\}\big]\le \E \big[\sup\{\sup_{s \in [0,t]} |\phi_{s}(x)-\phi_{s}(0)| \hspace{0.1cm} \colon  x \in \mathrm{B}(0,R)\}\big]\\&+\E \big[\sup_{s \in [0,t]} |\phi_{s}(0)|\big]
=\int\limits_0^\infty \P \big[\sup\{\sup_{s \in [0,t]}|\phi_{s}(x)-\phi_{s}(0)| \hspace{0.1cm} \colon  x \in \mathrm{B}(0,R)\}\ge z\big]\dd z+\E \big[\sup_{s \in [0,t]} |\phi_{s}(0)|\big].
\end{align*}
As $(\phi_{t}(x))$ is an Ornstein-Uhlenbeck process we certainly have a constant $L$ depending on $t$ such that $\E \big[\sup_{s \in [0,t]}|\phi_{s}(0)|\big] \le L$.\\

In the following we show the remaining part of Lemma~\ref{iouf_atr_lem3}, i.e. that there exists $K>0$, such that for all $R>0$
\begin{align*}
\int\limits_0^\infty \P\big[\sup\{\sup_{s \in [0,t]}|\phi_{s}(x)-\phi_{s}(0)| \hspace{0.1cm} \colon  x \in \mathrm{B}(0,R)\}> z\big] \dd z \le K R \,.
\end{align*}
Applying the chaining technique  to the function $ \sup\limits_{s\in [0,t]}|\phi_s(x)-\phi_s(0)| $ with suitable $\chi$,
\begin{align*}
\X=\mathrm{B}(0,R),\,\,\,\,\,\,\,0<\epsilon<1, \hspace{0.5cm} \epsilon_i=(1-\epsilon)\frac{e-1}{e}e^{-i} \hspace{0.5cm}\textrm{and}\hspace{0.5cm}\delta_i=R e^{-2i}
\end{align*}
we get for arbitrary $m>0$ and some $c_1>0$:
\begin{align*}
&\P \big[\sup\{\sup_{s \in [0,t]}|\phi_{s}(x)-\phi_{s}(0)| \hspace{0.1cm} \colon  x \in \mathrm{B}(0,R)\}> mR\big]\le \P \big[\sup_{s \in [0,t]}|\phi_{s}(0)-\phi_{s}(0)| > \epsilon m R \big]\\
&+  \sum_{i=0}^{\infty}|\chi_{i+1}|\sup_{|x-y|\le \delta_i}\P \big[\big|\sup_{s \in [0,t]}|\phi_{s}(x)-\phi_{s}(0)|-\sup_{s \in [0,t]}|\phi_{s}(y)-\phi_{s}(0)| \big|> \epsilon_i mR\big]\\
&\le c_1 \sum_{i=0}^{\infty}\left (\frac{R}{Re^{-2(i+1)}} \right )^d\sup_{|x-y|\le \delta_i}\P \big[\sup_{s \in [0,t]}\left\{|\phi_{s}(x)-\phi_{s}(0)|-|\phi_{s}(y)-\phi_{s}(0)|\right \}> \epsilon_i mR\big]\\
&\le c_1 \sum_{i=0}^{\infty}e^{2d(i+1)}\sup_{|x-y|\le \delta_i}\P\big[\sup_{s \in [0,t]}|\phi_{s}(x)-\phi_{s}(y)| > \epsilon_i mR\big]\,.
\end{align*}
According to Corollary \ref{iouf_gen_cor_0} there are constants $\lambda, \sigma>0$, such that  for some  standard Brownian motion $\left(B_t\right)_{t\ge 0}$ it holds almost surely  that
\begin{align*}
\sup_{s \in [0,t]}|\phi_{s}(x)-\phi_{s}(y)| \le |x-y|\exp(\sigma B_t^\star + \lambda t)\,.
\end{align*}
Therefore  we have
\begin{align*}
&\P\big [\sup\{\sup_{s \in [0,t]}|\phi_{s}(x)-\phi_{s}(0)| \hspace{0.1cm} \colon  x \in \mathrm{B}(0,R)\}> mR\big ]\\
&\le c_1 \sum_{i=0}^{\infty}e^{2d(i+1)}\P\big  [\delta_i \exp(\sigma B_t^\star + \lambda t) > \epsilon_i mR\big ]\le c_1 \sum_{i=0}^{\infty}e^{2d(i+1)}\P\big  [ B_1^\star  > \frac{1}{\sigma \sqrt t}\ln \frac{\epsilon_i mR}{\delta_i}-\frac{\lambda}{\sigma}\sqrt t\big ]\\
%&=c_1 \sum_{i=0}^{\infty}e^{2d(i+1)}\P\big  [ B_1^\star  > \frac{1}{\sigma \sqrt t}(\ln \frac{(1-\epsilon)(e-1) e^{-i}mR}{e R e^{-2i}}-\lambda t)\big  ]\\
&=c_1 \sum_{i=0}^{\infty}e^{2d(i+1)}\P\big  [ B_1^\star  > \frac{1}{\sigma \sqrt t}(\ln \frac{(1-\epsilon)(e-1) }{e }+\ln m+i-\lambda t)\big  ]\,.
\end{align*}
Let $m_0>0$ be so large that
\begin{align*}
\ln \frac{(1-\epsilon)(e-1) }{e }+\ln m_0-\lambda t>\sigma \sqrt t\,.
\end{align*}
Then for all $m>m_0$ with the help of Lemma \ref{comp_lemma_2} we get for some $c_2>0$
\begin{align*}
&\P\big  [\sup\{\sup_{s \in [0,t]}|\phi_{s}(x)-\phi_{s}(0)| \hspace{0.1cm} \colon  x \in \mathrm{B}(0,R)\}> mR\big ]\\[0.4cm]
&\le c_2 \sum_{i=0}^{\infty}\exp \big  ( 2d(i+1) -\frac{1}{2\sigma^2 t}(\ln \frac{(1-\epsilon)(e-1) }{e }+\ln m+i-\lambda t)^2 \big  )=:c_2\sum_{i=0}^\infty a_i\,,
\end{align*}
with
\begin{align*}
a_i=\exp \big  ( 2d(i+1) -\frac{1}{2\sigma^2 t}(\ln \frac{(1-\epsilon)(e-1) }{e }+\ln m + i -\lambda t)^2 \big  )\,.
\end{align*}
Choose $m_1(\lambda,\sigma,t,\epsilon)=m_1>m_0$ so big that for all $m>m_1$  the inequality
\begin{align*}
\frac{a_{i+1}}{a_i}=\exp \big  ( 2d -\frac{1}{\sigma^2 t}(\ln m + i+\frac{1}{2} -\lambda t+\ln \frac{(1-\epsilon)(e-1) }{e }) \big  ) \le e^{-1}
\end{align*}
holds
for all $i \in \N \cup \{0\}$ (that is, the inequality holds  for $i=0$).
Then for all $m>m_1$ we have
\begin{align*}
&\P\big [\sup\{\sup_{s \in [0,t]}|\phi_{s}(x)-\phi_{s}(0)| \hspace{0.1cm} \colon  x \in \mathrm{B}(0,R)\}> mR\big ]\\[0.4cm]
& \le c_2 \sum_{i=0}^{\infty}a_i \le c_2 a_0 \sum_{i=0}^{\infty}e^{-i} = c_2\frac{e}{e-1} \exp \big  ( 2d -\frac{1}{2\sigma^2 t}(\ln \frac{(1-\epsilon)(e-1) }{e }+\ln m  -\lambda t)^2 \big  )
\end{align*}
and clearly there are positive constants $p(t,\sigma,\lambda,d)$ and $q(t,\sigma,\lambda,d)$,  such that for all $m\ge m_1$
\begin{align}\label{cr_short}
&\P\big [\sup\{\sup_{s \in [0,t]}|\phi_{s}(x)-\phi_{s}(0)| \hspace{0.1cm} \colon  x \in \mathrm{B}(0,R)\}> mR\big ] \le p \exp \big  ( -q\ln^2 m  \big  ) \,.
\end{align}
Using the last inequality we finally get
\begin{align*}
&\E \big [\sup\{\sup_{s \in [0,t]}|\phi_{s}(x)-\phi_{s}(0)| \hspace{0.1cm} \colon  x \in \mathrm{B}(0,R)\}\big ]\\
&=R\int\limits_0^\infty \P\big [\sup\{\sup_{s \in [0,t]}|\phi_{s}(x)-\phi_{s}(0)| \hspace{0.1cm} \colon  x \in \mathrm{B}(0,R)\}\ge mR\big ]\dd m\\
&\le R m_1 + R\int\limits_{m_1}^\infty \P \big [\sup\{\sup_{s \in [0,t]}|\phi_{s}(x)-\phi_{s}(0)| \hspace{0.1cm} \colon  x \in \mathrm{B}(0,R)\}\ge mR\big ]\dd m\\
&\le R m_1 + R\int\limits_{m_1}^\infty p \exp \left ( -q\ln^2 m  \right ) \dd m \le KR
\end{align*}
for some $K>0$.
\end{proof}
The following corollary is a consequence of the preceding two lemmata
\begin{corollary}$ $
\begin{itemize}
\item[\textnormal{(i)}] For every $t>0$ there exists $R_0>0$ (depending on $t$),  such that for all  $R \ge R_0$ and all $x\in \R^d$
\begin{align}\label{diam_estimate}
\P \big(\textnormal{diam}(\phi_t(\textnormal{B}(x,1)))\ge  R\big )\le c_1(t,\sigma,\lambda,d) \exp \left ( -c_2(t,\sigma,\lambda,d)\ln^2 R\right)\,
\end{align}
holds for all $R\ge R_0$ with some  positive constants $c_1(t,\sigma,\lambda,d)$ and $c_2(t,\sigma,\lambda,d)$ (depending on $t$, $\sigma$, $\lambda$ and $d$).
\item[\textnormal{(ii)}] For all $\gamma_0>0$ there exist $t_0>0$ and $R_0>0$ such that for all $\gamma\ge \gamma_0$ and $R \ge R_0$
\begin{align}\label{iouf_atr_prop1}
\P\big(\sup\{ | \phi_{t_0}(x) | \hspace{0.1cm} \colon  x \in \textnormal{B}(0,R) \} >\gamma R)\le c_1  \exp \big ( -c_2\ln^2(\gamma R)\big ) ,
\end{align}
where  $c_1$ and $c_2$ are positive constants depending on $d$, $\lambda$, $\sigma$, $c$ and $t_0$.
\end{itemize}
\end{corollary}
\begin{proof}
(i) follows immediately from (\ref{cr_short})  (ii) is a combination of (i) and Lemma~\ref{O_U_estimate}.
\end{proof}

\begin{propo}\label{iouf_atr_prop2}
For all $\delta \in (0,1)$ there exist $t_0>0$ and $R_0>0$, such that for all  $R>0$
\begin{align*}
\E\left[\sup\{|\phi_{t_0}(x)| \hspace{0.1cm} \colon  x\in \textnormal{B}(0,R)\}\right] \le \delta R \eins_{\{R \ge R_0\}}+\delta R_0 \eins_{\{R < R_0\}}\,.
\end{align*}
\end{propo}
\begin{remark}
 For our purposes the existence of one such  $\delta$ is enough.
\end{remark}
\begin{proof}

\textbf{Step 1:} Choose arbitrary  $\delta \in (0,1)$. We first show that there exist  $t_0>0$ and $R_0>0$,  such that for all $R \ge R_0$
\begin{align*}
\E\left[\sup\{|\phi_{t_0}(x)| \hspace{0.1cm} \colon  x\in \mathrm{B}(0,R)\}\right] \le \delta R\,.
\end{align*}
Let  $\gamma_0=\frac{\delta}{2}$. Choose $R_{0,1}>0$ and $t_0>0$ according to  \eqref{iouf_atr_prop1}. Thus  for all $R \ge R_{0,1}$
\begin{align*}
&\E\left[ \sup\{|\phi_{t_0}(x)| : x\in \mathrm{B}(0,R)\} \right]= R \int\limits_0^\infty \P(\sup\{|\phi_{t_0}(x)| \hspace{0.1cm} \colon  x\in \mathrm{B}(0,R)\}>\gamma R)\,\dd\gamma\\
&\le R \gamma_0+R \int\limits_{\gamma_0}^ \infty \P(\sup\{|\phi_{t_0}(x)| \hspace{0.1cm} \colon  x\in \mathrm{B}(0,R)\}>\gamma R)\,\dd\gamma \le R \gamma_0+R\int\limits_{\gamma_0}^ \infty c_1   \exp \big ( -c_2\ln^2(\gamma R)\big ) \,\dd\gamma,
\end{align*}
where  $c_1$ and $c_2$ are as in  \eqref{iouf_atr_prop1}. Choose  $R_{0,2}>0$   such that
\begin{align*}
\int\limits_{\gamma_0}^\infty  c_1 \exp \big ( -c_2\ln^2 (\gamma R)\big )\,\dd s \le \gamma_0,
\end{align*}
which is possible as the whole expression on the left-hand side goes to zero with $R\to +\infty$.  Clearly,  for all $R>R_0:=\max\{R_{0,1},R_{0,2}\}$ it holds that $\E\big[\sup\{|\phi_{t_0}(x) :   x\in \mathrm{B}(0,R)\}\big] \le  \delta R$.\\
\textbf{Step 2:} For arbitrary $R>0$, using the result of Step 1 we have
\begin{align*}
&\E\left[\sup\{|\phi_{t_0}(x) :   x\in \mathrm{B}(0,R)\}\right] \\
 &= \eins_{\{R \ge R_0\}} \E\left[\sup\{|\phi_{t_0}(x) :  x\in \mathrm{B}(0,R)\}\right]+\eins_{\{R < R_0\}}\E\left[\sup\{|\phi_{t_0}(x) :  x\in \mathrm{B}(0,R)\}\right]\\
&\le \eins_{\{R \ge R_0\}}\delta R + \eins_{\{R < R_0\}}\E\left[\sup\{|\phi_{t_0}(x)| \hspace{0.1cm} \colon  x\in \mathrm{B}(0,R_0)\}\right]\le \eins_{\{R \ge R_0\}}\delta R + \eins_{\{R < R_0\}}\delta R_0,
\end{align*}
where $R_0$ is of course as in Step 1.
\end{proof}
Now we give the proof of the main result:\\
\textbf{Proof of the main theorem.}\\
In the following we verify the sufficient condition for the existence of a random weak attractor stated in Corollary \ref{iouf_atr_cor}, i.e. for an IOUF $\phi$ we show that  there exists a constant $M>0$, such that for every $R>0$ there exists a  $T>0$ with
\begin{align*}
\E\big[\sup\{|\phi_{t}(x)|  : x\in \mathrm{B}(0,R)\}\big] \le M \hspace{0.7cm}\textrm{ for all }\hspace{0.5cm} t\ge T.
\end{align*}
We do this in two steps:\\
\textbf{Step 1:} Choose arbitrary $\delta \in (0,1)$. Take $t_0>0$ and $R_0>0$ as in Proposition \ref{iouf_atr_prop2}, that is for all $R>0$,
$$\E\big[\sup\{|\phi_{t_0}(x)| : x\in \mathrm{B}(0,R) \}\big]\le \delta R \eins_{\{R\ge R_0\}}+\delta R_0\eins_{\{R< R_0\}}\,. $$
 We first show that there exists a  constant $M_1>0$, such that for every $R>0$ there is an $n_0 \in \N$ such that for all $n \ge n_0$
\begin{align*}
\E\big[\sup\{|\phi_{n t_0}(x)| \hspace{0.1cm} \colon  x\in \mathrm{B}(0,R)\}\big] \le M_1 \textrm{ for all } n\in \N.
\end{align*}
Observe that for arbitrary $n\ge 1$
\begin{align*}
&\E\big [\sup\{|\phi_{n t_0}(x)| \hspace{0.1cm} \colon  x\in \mathrm{B}(0,R)\}\,\big | \,\F_{(n-1)t_0} \big ]\\[0.4cm]
& =\E\big [\sup\{|\phi_{(n-1)t_0,n t_0}(x)| \hspace{0.1cm} \colon  x\in \phi_{(n-1)t_0}(\mathrm{B}(0,R))\}\,\big |\, \F_{(n-1)t_0} \big]=(\star)
\end{align*}
and since $\phi_{(n-1)t_0}$ is $\F_{(n-1)t_0}$-measurable and $\phi_{(n-1)t_0,n t_0}$ is independent of $\F_{(n-1)t_0}$ we obtain
\begin{align*}
(\star)& =\E\big [\sup\{|\phi_{(n-1)t_0,n t_0}(x)|\hspace{0.1cm} \colon  x\in B \} \big ] \Big |_{B=\phi_{(n-1)t_0}(\mathrm{B}(0,R))}\\[0.4cm]
& \le \E\big [\sup\{|\phi_{(n-1)t_0,n t_0}(x)| \hspace{0.1cm} \colon  x\in \mathrm{B}(0,\tilde R)\} \big ]\Big |_{\tilde R=\sup\{|\phi_{(n-1)t_0}(x)| \hspace{0.1cm} \colon  x \in \mathrm{B}(0,R))\}}\\[0.4cm]
&=\E\big [\sup\{|\phi_{t_0}(x)| \hspace{0.1cm} \colon  x\in \mathrm{B}(0,\tilde R)\} \big ]\Big |_{\tilde R=\sup\{|\phi_{(n-1)t_0}(x)| \hspace{0.1cm} \colon  x \in \mathrm{B}(0,R))\}}
\end{align*}
where the last equality is due to the  time homogeneity of the flow. Applying  Proposition~\ref{iouf_atr_prop2} we get
\begin{align*}
&\E\big [\sup\{|\phi_{n t_0}(x)| \hspace{0.1cm} \colon  x\in \mathrm{B}(0,R)\}\,\big|\, \F_{(n-1)t_0} \big ]\\
& \le \E\big [\sup\{|\phi_{t_0}(x)| \hspace{0.1cm} \colon  x\in \mathrm{B}(0,\tilde R)\} \big ]\Big |_{\tilde R=\sup\{|\phi_{(n-1)t_0}(x)| \hspace{0.1cm} \colon  x \in \mathrm{B}(0,R))\}}\\
& \le \Big [  \delta \tilde R \eins_{\{\tilde R\ge R_0\}}+ \delta  R_0 \eins_{\{\tilde R< R_0\}}\Big ]_{\tilde R=\sup\{|\phi_{(n-1)t_0}(x)| \colon  x \in \mathrm{B}(0,R))\}}%\\[0.4cm]
%&=  \delta \sup\{|\phi_{(n-1)t_0}(x)| : x \in \mathrm{B}(0,R))\}  \eins_{\{\sup\{|\phi_{(n-1)t_0}(x)| \hspace{0.1cm} \colon  x \in \mathrm{B}(0,R))\} \ge R_0\}}\\[0.4cm]
%&\hspace{0.4cm}+\delta R_0\eins_{\{\sup\{|\phi_{(n-1)t_0}(x)| \hspace{0.1cm} \colon  x \in \mathrm{B}(0,R))\} < R_0\}}\,.
\end{align*}
That is, with $X_n^R:=\sup\{|\phi_{nt_0}(x)| : x \in \mathrm{B}(0,R))\}$ and $\G_n:=\F_{n t_0}$ we have
\begin{align}\label{iouf_atr_prop3.1}
\E\big[X_n^R\,\big |\,\G_{n-1}\big]\le \delta X_{n-1}^R\eins_{\{X_{n-1}^R\ge R_0\}}+\delta R_0 \eins_{\{X_{n-1}^R< R_0\}}\,.
\end{align}
Now   the following iteration argument:
\begin{align*}
&\E[X_n^R]=\E[ \dots \E[\E[X_n^R\,\big|\,\G_{n-1}]|\G_{n-2}] \dots \big |\,\G_0] \\
&\le \E[ \dots \E[\delta X_{n-1}^R\eins_{\{X_{n-1}^R\ge R_0\}}+\delta R_0 \eins_{\{X_{n-1}^R< R_0\}}\,\big|\,\G_{n-2}] \dots \big|\,\G_0] \\
&\le \delta\E[ \dots \E[X_{n-1}^R\,\big|\,\G_{n-2}] \dots \big|\,\G_0]+ \delta R_0 \le \delta^2\E[ \dots \E[X_{n-2}^R\,\big|\,\G_{n-3}] \dots \big|\,\G_0]+ \delta^2 R_0 +\delta R_0 \le \dots\\
& \le \delta^{n-1}\E X_1^R+R_0\delta (1+\delta + \dots +\delta^{n-1}) \\
&\le \delta^n R \eins_{\{R \ge R_0\}}+ \delta^n R_0 \eins_{\{R < R_0\}}+R_0 \frac{\delta}{1-\delta}\le \delta^n \max \{R,R_0 \}+R_0 \frac{\delta}{1-\delta}
\end{align*}
 proves the statement, since obviously $\limsup_{n \to \infty}\E \big[ X_n^R \big]\le \frac{R_0 \delta}{1-\delta}$, and  this estimate is independent of $R$.\\
\textbf{Step 2:} Let $t=nt_0+s$ with $s<t_0$ for $t_0$ as in Step 1. \\
An analogous argument as in Step 1  together with Lemma \ref{iouf_atr_lem3} implies:
\begin{align*}
&\E\big[\sup\{|\phi_{t}(x)| \hspace{0.1cm} \colon  x\in \mathrm{B}(0,R)\}\big] = \E\left [ \E\big[\sup\{|\phi_{t}(x)| \hspace{0.1cm} \colon  x\in \mathrm{B}(0,R)\}\,\big |\, \F_{nt_0}\big]\right ]\\[0.3cm]
&=\E\left [ \E\big[\sup\{|\phi_{nt_0,t}(x)| \hspace{0.1cm} \colon  x\in \phi_{nt_0}(\mathrm{B}(0,R))\}\,\big|\, \F_{nt_0}\big]\right ]\\[0.3cm]
&\le\E\left [  \E\big[\sup\{|\phi_{s}(x)| : x\in \mathrm{B}(0,\tilde R)\}\big]\Big |_{\tilde R=\sup\{|\phi_{nt_0}(x)| \hspace{0.1cm} \colon  x\in \mathrm{B}(0,R)\}}\right ]\\[0.3cm]
&\le\E\left [  \E\big[\sup\{\sup_{0\le u\le t_0}|\phi_{u}(x)| : x\in \mathrm{B}(0,\tilde R)\}\big]\Big |_{\tilde R=\sup\{|\phi_{nt_0}(x)| \hspace{0.1cm} \colon  x\in \mathrm{B}(0,R)\}}\right ]\\[0.3cm]
&\le\E \left [ K(t_0)\sup\{ |\phi_{nt_0}(x) | \hspace{0.1cm} \colon  x\in \mathrm{B}(0,R)\} +L(t_0) \right ]  \hspace{3cm} (\text{Lemma \ref{iouf_atr_lem3}})\\[0.3cm]
&=  K(t_0) \E \big[\sup\{ |\phi_{nt_0}(x) | \hspace{0.1cm} \colon  x\in \mathrm{B}(0,R)\}\big] +L(t_0)\\[0.3cm]
& \le K(t_0) \delta^n \max\{R,R_0\}+K(t_0) R_0 \frac{\delta}{1-\delta}  +L(t_0)\hspace{4.2cm} (\text{Step 1})
\end{align*}
and therefore for all $R>0$
\begin{align*}
&\limsup_{t \to +\infty}\E\big[\sup\{|\phi_{t}(x)| \hspace{0.1cm} \colon  x\in \mathrm{B}(0,R)\}\big]  \le K(t_0) R_0 \frac{\delta}{1-\delta}  +L(t_0).
\end{align*}
The statement is now proven since the  bound above  is independent of $ R $. \hfill $\square $

\subsection{The case of negative Lyapunov exponent}
 Our aim here is to show that in the case of negative top Lyapunov exponent the global weak set attractor of an isotropic Ornstein-Uhlenbeck flow   is trivial, i.e. contains almost surely only one point. We sketch the proof under the following additional  ''squeezing`` condition.\\
\begin{condition}(Squeezing)\\
We say a stochstic flow $\phi$ satisfies the squeezing condition if  for every $r>0$ there exists $t_r>0$, such that the ball with radius $r$ and centered at zero is uniformly squeezed by $\phi$ in the sense that  there exists $\epsilon>0$ with the property that
\begin{align*}
\P(\phi_{0,t_r}{\text{Ball}(0,r+\epsilon)}\subset\text{Ball}(0,r-\epsilon) )>0.
\end{align*}
\end{condition}

The squeezing condition can be verified for  IBFs whose potential measure in the spectral decomposition of the corresponding field $F$ is not supported on the set of zeros of certain Bessel functions. %For details consult \cite{Dim06}, Theorems 6.2.1 and 6.2.2.
Moreover, an ongoing joint work of one of the authors  with Steffen Dereich on a support theorem for stochastic flows shows that the squeezing condition is always verified by an IOUF.\\

\begin{theo}\label{iouf_atr_triv}
Let $\phi$ be an isotropic Ornstein-Uhlenbeck flow with generating field\\ $F(t,x)-cxt$ and strictly negative top Lyapunov exponent verifying the squeezing condition.  Then the weak attractor $A(\omega)$ of $\phi$ is trivial, i.e.  contains almost surely only one point.
\end{theo}
\begin{proof}
We will divide the proof into   several steps.\\
\textbf{Step 1}: Simple iteration of the squeezing condition and the time homogeneity of the IOUF imply that for all $R>r>0$ there exists $T_\phi^{R,r}>0$, such that
\begin{align}\label{iouf_atr_triv0}
\P(\phi_{nT_\phi^{R,r},(n+1)T_\phi^{R,r}}{\text{Ball}(0,R)}\subset\text{Ball}(0,r) )>0 \,\,\,\,\text{ for all } n\ge 0
\end{align}
\textbf{Step 2}: Here we pass to a certain two-sided RDS $\psi$, corresponding to the isotropic Ornstein-Uhlenbeck flow seen from the point of view of a moving particle. \\
Consider the semimartingale field $\left(W(t,x)\colon t\in \R_+, x\in \R^d \right)$ defined as the $C^{3,\delta}$-valued modification (for arbitrary $\delta<1$) of the field
\begin{align*}
\int \limits_0^t V(\dd u,x+\phi_{0,u}(0))-\int\limits_0^tV(\dd u,\phi_{0,u}(0))=-c\int\limits_0^t x \dd u + \int\limits_0^t F(\dd u, x+\phi_{0,u}(0))- \int\limits_0^t F(\dd u, \phi_{0,u}(0))\,.
\end{align*}
The  joint variation process is given by $\big(b_{i,j}(x-y)+b_{i,j}(0)-b_{i,j}(x)-b_{i,j}(y)\big) t$
and since it is  deterministic  the martingale part of
$\left\{W(t,x)\colon t\in \R_+, \,x\in \R^d \right\}$ is a Gaussian martingale field. Moreover, $W$ has stationary and independent increments. Let $\psi$ be  the forward Brownian  stochastic flow  generated by $W(x,t)$ via the SDE:
\begin{align}\label{iouf_watr_triv_thm1}
\psi_{s,t}(x)-x=\int\limits_s^t W(\dd u,\psi_{s,u}(x)) \,\,\,\, \text{ for all } 0\le s \le t\,\,\text{ and } x\in \R^d\,.
\end{align}
With the help of Theorem 3.3.3 in \cite{Ku90} and the pathwise uniqueness of the solutions of the SDE above we can show that actually
\begin{align}\label{iouf_watr_triv_thm2}
\psi_{s,t}(x)=\phi_{s,t}(x+\phi_{0,s}(0))-\phi_{0,t}(0)\,,
\end{align}
i.e. $\psi$ is simply the flow $\phi$, observed from the point of view of the moving particle started at the origin. Of course  the last equality is in the sense that the flows on the left and the right hand side are modifications of each other. \\
As in~\cite{AS95} we can construct a filtered probability space $(\tilde \Omega, \tilde \F,  ( \tilde\F_s^t  )_{-\infty< s\le t < +\infty}, \tilde \P)$, an ergodic filtered metric dynamical system $\left(\theta_t\right )_{ t \in \R}$ and a perfect cocycle $\psi\colon  \R \times \R^d  \times \tilde\Omega \to \R^d$ (we use the same symbol $\psi$ for the RDS and the flow!), such that the distribution of $\left \{ \psi(t-s,\theta_s \tilde \omega) : s,t \ge 0\right \}$ coincides with the distribution of $\left \{\psi_{s,t} :   s, t\ge 0 \right \}$.\\
It is also easily concluded   that for arbitrary $R>r>0$ there is a $T^{R,r}>0$, such that with positive probability the ball with Radius $R$ will be squeezed after time $T^{R,r}$ by the action of the cocycle $\psi$ into the ball with radius $r$, i.e.
\begin{align}\label{iouf_atr_triv1}
\tilde \P\left (\psi(T^{R,r},\tilde \omega)(\mathrm{B}(0,R))\subset \mathrm{B}(0,r)\right )>0  \,.
\end{align}
\textbf{Step 3}:
Here we sketch the proof of the fact  that the global stable manifold $\mathcal{S}_g$ of the RDS $\psi$ contains with positive probability a ball centered at zero and having arbitrarily large radius. The idea  is simple: the distance $(d_t)$ between the stationary (here fixed) trajectory and the complement of the global stable manifold (see \cite{MoSch99}) is stationary. However, according to Step 1, the image of arbitrary big ball becomes arbitrary small with strictly positive probability in finite time $t_0$, that is with positive probability at $t_0$ the image of the ball is contained in the global stable manifold, hence the stationary process $(d_t)$ assumes arbitrarily large values with strict positive probability.  \\
\textbf{Step 4}:
The RDS $\psi$ has a weak attractor, since the bound in Corollary \ref{iouf_atr_cor} holds also for $\psi$ if it holds for $\phi$. Denote by $A(\tilde \omega)$ the attractor of the RDS $\psi$. Similar argument as above implies that $A(\tilde \omega)$ is contained in the global stable manifold and thus is trivial. The same must hold for the random attractor of $\phi$ as it is a random shift of the attractor of $\psi$ (see~\eqref{iouf_watr_triv_thm2}) and  clearly a random translation of a set containing only one random point is still a random one-point set.     .
\end{proof}

\bibliography{biblio}
\bibliographystyle{plain}
\bibliographystyle{abbrv}

\end{document}